\documentclass[fleqn,a4paper,10pt]{amsart}

\usepackage{graphicx}
\usepackage{tikz}
\usetikzlibrary{shapes,decorations.pathmorphing}
\usepackage{amsmath}
\usepackage{amssymb,amsthm}
\usepackage{comment}
\usepackage{hyperref}
\usepackage{enumerate}

\theoremstyle{definition}

\title[On a conjecture by Chapuy about Vorono\"\i\ cells in large maps]{On a conjecture by Chapuy about Vorono\"\i\ cells in large maps}
\author{Emmanuel Guitter}
\address{Institut de physique th\'eorique, Universit\'e Paris Saclay, CEA, CNRS, F-91191 Gif-sur-Yvette}
\email{emmanuel.guitter@cea.fr}

\begin{document}
\maketitle

\begin{abstract}
In a recent paper, Chapuy conjectured that, for any positive integer $k$, the law for the fractions of total area covered by the $k$ Vorono\"\i\ cells 
defined by $k$ points picked uniformly at random in the Brownian map of any fixed genus is the same law as that of a uniform $k$-division of
the unit interval. For $k=2$, i.e.\ with two points chosen uniformly at random, it means that the law for the ratio of the area of one of the two Vorono\"\i\ cells
by the total area of the map is uniform between $0$ and $1$. Here, by a direct computation of the desired law, we show that this latter conjecture for $k=2$ 
actually holds in the case of large planar (genus $0$) quadrangulations as well as for large general planar maps (i.e.\ maps whose faces have
arbitrary degrees). This corroborates Chapuy's conjecture in its simplest realizations.
 \end{abstract}

\section{Introduction}
\label{sec:introduction}

The asymptotics of the number of maps of some arbitrary given genus has been known for quite a while \cite{BC86} and involves some universal constants
$t_g$, whose value may be determined recursively. In its simplest form, the $t_g$-recurrence is a simple quadratic recursion for the $t_g$'s,
first established in the physics literature \cite{GM90,DS90,BK90} in the context of matrix integrals, then proven rigorously in the mathematical literature \cite{BGR08,GJ08,CC15}.
In a recent paper \cite{GC16}, Chapuy addressed the question of reproducing the $t_g$-recurrence in a purely combinatorial way. 
By a series of clever arguments involving various bijections, he could from his analysis extract exact values for a number of moments of the law for
the area of Vorono\"\i\ cells defined by uniform points in the Brownian map of some arbitrary fixed genus.
In view of these results and other evidence, he was eventually led to formulate the following conjecture: \emph{for any integer $k\geq 2$, the proportions of the total area 
covered by the $k$ Vorono\"\i\ cells defined by $k$ points picked uniformly at random in the Brownian map of any fixed genus have the same law as a uniform $k$-division of
the unit interval}. The simplest instance of this conjecture is for the planar case (genus $0$) and for $k=2$. It may be rephrased by 
saying that, given two points picked uniformly at random in the planar Brownian map and the corresponding two Vorono\"\i\ cells, \emph{the law for the ratio of the area of one of the 
cells by the total area of the map is uniform between $0$ and $1$}. 

The aim of this paper is to show that this latter conjecture ($k=2$ and genus $0$) is actually true by computing the desired law for particular realizations of the planar Brownian map, 
namely large random planar quadrangulations and large random general planar maps (i.e.\ maps whose faces have
arbitrary degrees). We will indeed show that, for planar quadrangulations with a fixed area $N$ ($=$ number of faces) and with two marked vertices 
picked uniformly at random, the law for ratio $\phi=n/N$ between the area $n$ of the Vorono\"\i\ cell around, say, the second vertex and the total area $N$ is, for large $N$ 
and finite $\phi$, the uniform law in the interval $[0,1]$. This property is derived by a \emph{direct computation of the law} itself from explicit discrete or asymptotic enumeration results.
The result is then trivially extended to Vorono\"\i\ cells in general planar maps of large area (measured in this case by the number of edges).

\section{Vorono\"\i\ cells in bi-pointed quadrangulations}
\label{sec:voronoi}
\begin{figure}
\begin{center}
\includegraphics[width=8cm]{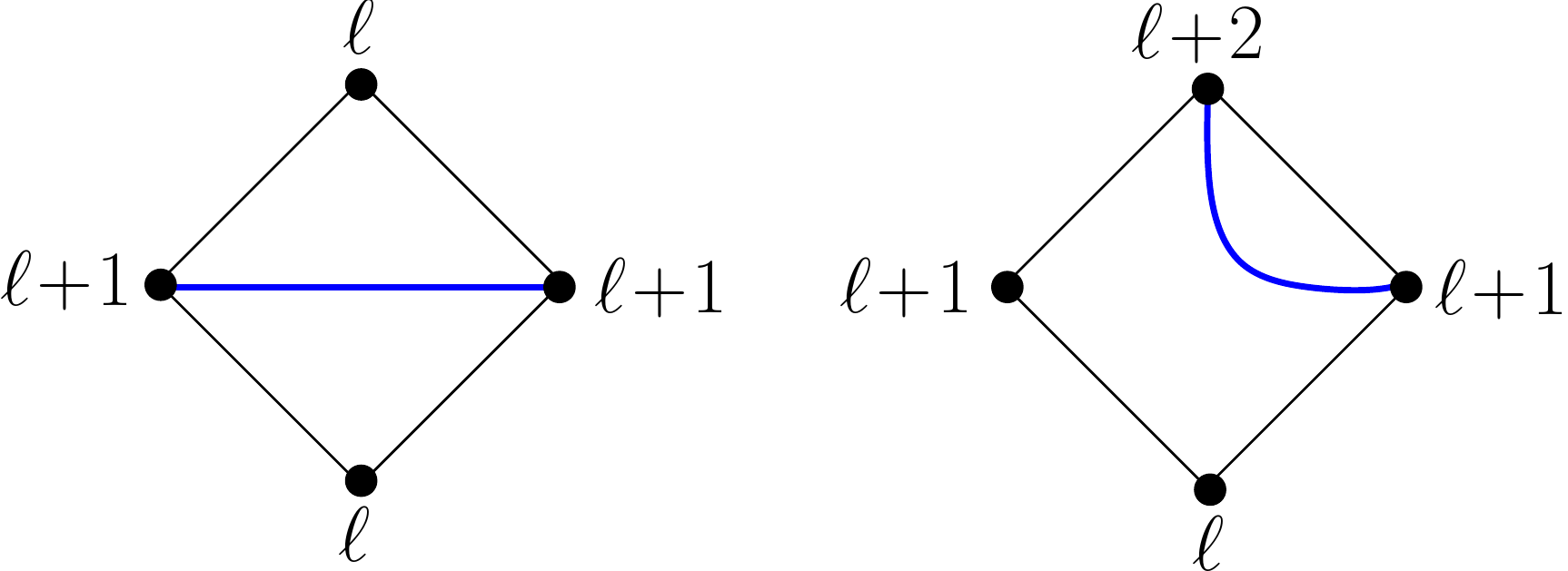}
\end{center}
\caption{The local rules of the Miermont bijection. These rules are the same as those of the Schaeffer bijection.}
\label{fig:SR}
\end{figure}
\begin{figure}
\begin{center}
\includegraphics[width=8cm]{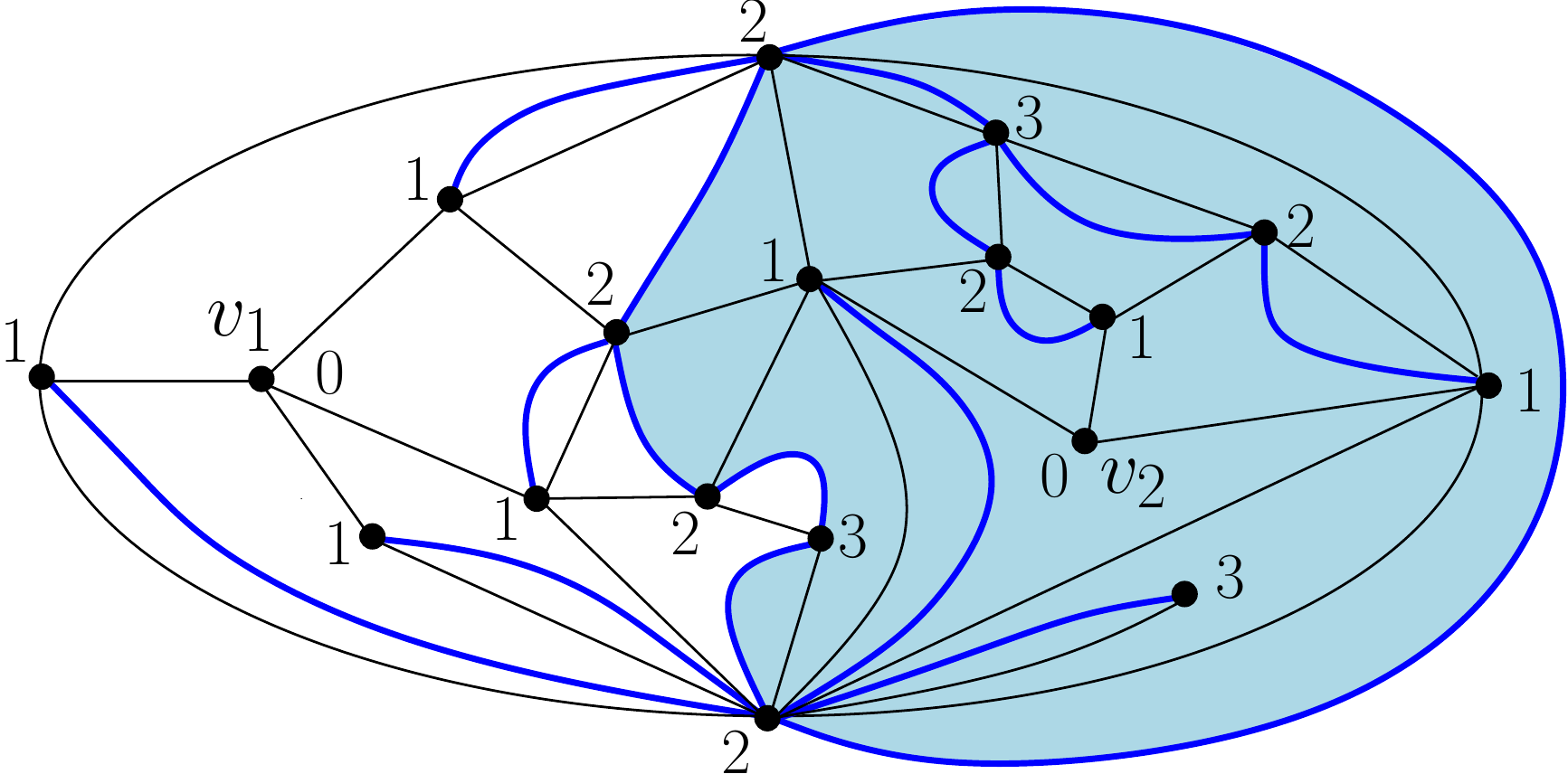}
\end{center}
\caption{A bi-pointed planar quadrangulation and the associated i-l.2.f.m via the Miermont bijection. The two faces of the i-l.2.f.m delimit
two domains on the quadrangulation which define our two Vorono\"\i\ cells. Here, one of the cells has been filled in light-blue.}
\label{fig:cells}
\end{figure}
This paper deals exclusively with planar maps, which are connected graphs embedded on the sphere. Our starting point are \emph{bi-pointed planar quadrangulations}, which
are planar maps whose all faces have degree $4$, and with two marked distinct vertices, distinguished as $v_1$ and $v_2$. For convenience, we will assume 
here and throughout the paper that \emph{the graph distance $d(v_1,v_2)$ between $v_1$ and $v_2$ is even}. As discussed at the end of Section~\ref{sec:scaling}, this requirement is not crucial
but it will make our discussion slightly simpler. The Vorono\"\i\ cells associated to $v_1$ and $v_2$ regroup, so to say, the set of vertices which are closer to 
one vertex than to the other. A precise definition of the Vorono\"\i\ cells in bi-pointed planar quadrangulations may be given upon coding these maps via the
well-known Miermont bijection\footnote{We use here a particular instance of the Miermont bijection for two ``sources" and with vanishing ``delays".} \cite{Miermont2009}. It goes as follows: we first assign to each vertex $v$ of the quadrangulation
its label $\ell(v)=\min(d(v,v_1),d(v,v_2))$ where $d(v,v')$ denotes the graph distance between two vertices $v$ and $v'$ in the quadrangulation. The label $\ell(v)$ is thus the distance from $v$ to 
the closest marked vertex $v_1$ or $v_2$. The labels are non-negative integers
which satisfy $\ell(v)-\ell(v')=\pm 1$ if $v$ and $v'$ are adjacent vertices. Indeed, it is clear from their definition that labels between adjacent vertices can differ by at most $1$.
Moreover, a planar quadrangulation is bipartite so we may color its vertices in black and white in such a way that adjacent vertices carry different colors. Then if we chose
$v_1$ black, $v_2$ will also be black since $d(v_1,v_2)$ is even. Both $d(v,v_1)$ and $d(v,v_2)$ are then simultaneously even if $v$ is black and so is thus $\ell(v)$. 
Similarly, $d(v,v_1)$, $d(v,v_2)$ and thus $\ell(v)$ are odd if $v$ is white so that the parity of labels changes between adjacent neighbors. We conclude that labels between adjacent vertices necessarily differ by $\pm 1$.

The cyclic sequence of labels around a face is then necessarily of one of the two types displayed in Figure~\ref{fig:SR}, namely, if $\ell$ is the smallest label around the face, of
the form $\ell\to\ell+1\to\ell\to\ell+1$ or $\ell\to\ell+1\to\ell+2\to\ell+1$. Miermont's coding is similar to that of the well-known Schaeffer bijection \cite{SchPhD} and consists in drawing inside each face an edge connecting the two corners within
the face which are followed \emph{clockwise} by a corner with smaller label (here the label of a corner is that of the incident vertex). Removing all the original edges, we obtain a graph embedded 
on the sphere whose vertices are {\it de facto} labelled by integers (see Figure~\ref{fig:cells}). It was shown by Miermont \cite{Miermont2009} that this graph spans all the original vertices of the quadrangulation but $v_1$ and $v_2$, is connected and defines a planar map \emph{with exactly $2$ faces} $f_1$ and $f_2$, where $v_1$ (which is not part of the two-face map) lies strictly inside $f_1$, and $v_2$ strictly inside $f_2$. As for the 
vertex labels on this two-face map, they are easily shown to satisfy:
\begin{enumerate}[$\langle \hbox{a}_1\rangle$]
\item{Labels on adjacent vertices differ by $0$ or $\pm1$.}
\item{The minimum label for the set of vertices incident to $f_1$ is $1$.}
\item{The minimum label for the set of vertices incident to $f_2$ is $1$.}
\end{enumerate} 
In view of this result, we define a planar \emph{iso-labelled two-face map} (i-l.2.f.m) as a planar map with exactly two faces, \emph{distinguished} as $f_1$ and $f_2$, and whose vertices carry integer labels
satisfying the constraints $\langle \hbox{a}_1\rangle$-$\langle \hbox{a}_3\rangle$ above. Miermont's result is that the construction presented above actually provides a \emph{bijection between bi-pointed planar quadrangulations 
whose two distinct and distinguished marked vertices are at some even graph distance
from each other and planar i-l.2.f.m}. 
Moreover, the Miermont bijection guarantees that (identifying the vertices $v$ of the i-l.2.f.m with their pre-image in the associated quadrangulation):
\begin{itemize}
\item{The label $\ell(v)$ of a vertex $v$ in an i-l.2.f.m corresponds to the minimum distance $\min(d(v,v_1),d(v,v_2))$ from $v$ to the marked vertices 
$v_1$ and $v_2$ in the associated bi-pointed quadrangulation.} 
\item{All the vertices incident to the first face $f_1$ (respectively the second face $f_2$) in the i-l.2.f.m are closer to $v_1$ than to $v_2$ (respectively closer to $v_2$ than to $v_1$) 
or at the same distance from both vertices in the associated quadrangulation.}
\item{The minimum label $s$ among vertices incident to both $f_1$ and $f_2$ and the distance $d(v_1,v_2)$ between the marked vertices in the associated quadrangulation
are related by $d(v_1,v_2)=2s$.}
\end{itemize} 
Clearly, all vertices incident to both $f_1$ and $f_2$ are at the same distance from both $v_1$ and $v_2$. Note however that the reverse is not true and that 
vertices at equal distance from both $v_1$ and $v_2$ might very well lie strictly inside a given face. 

Nevertheless, the coding of bi-pointed quadrangulations by i-l.2.f.m provides us
with a well defined notion of Vorono\"\i\ cells. Indeed, since it has exactly two faces, the i-l.2.f.m is made of a simple closed loop separating the two faces, completed by 
(possibly empty) subtrees 
attached on both sides of each of the loop vertices (see Figure~\ref{fig:cells}).  Drawing the quadrangulation and its associated i-l.2.f.m on the same picture,
we \emph{define} the two Vorono\"\i\ cells of a bi-pointed quadrangulation
as the two domains obtained by cutting along the loop of the associated i-l.2.f.m. Clearly, each Vorono\"\i\ cell contains only vertices closer from one of the marked vertices that  
from the other (or possibly at the same distance). As just mentioned, vertices at the border between the two cells are necessarily at the same distance from $v_1$ and $v_2$. 
Note also that all the edges of the quadrangulation lie strictly in one cell or the other. This is not the case for all the faces of the quadrangulation whose situation
is slightly more subtle. Clearly, these faces are in bijection with the edges of the i-l.2.f.m. The latter come in three species, those lying strictly inside the first face of the i-l.2.f.m, in which case the associated face in the quadrangulation lies strictly inside the first cell, those lying strictly inside the second face of the i-l.2.f.m, in which case the associated face in the quadrangulation lies strictly inside the second cell, and those belonging to the loop separating the two faces of the i-l.2.f.m, in which case the associated face in the quadrangulation is split in two by the cutting and shared 
by the two cells.

If we now want to measure the area of the Vorono\"\i\ cells, i.e.\ the number of faces which they contain, several prescriptions may be taken to properly account for the shared faces.
The simplest one is to count them as half-faces, hence contributing a factor $1/2$ to the total area of each of the cells. For generating functions, this prescription amounts to assign a weight 
$g$ per face strictly within the first Vorono\"\i\ cell, a weight $h$ per face strictly within the second cell and a weight $\sqrt{g\, h}$ per face shared by the two cells. 
A different prescription would consist in attributing 
each shared face to one cell or the other randomly with probability $1/2$ and averaging over all possible such attributions. In terms of generating functions,
this would amount to now give a weight $(g+h)/2$ to the faces shared by the two cells. As discussed below, the precise prescription for shared faces turns out to be irrelevant 
in the limit of large quadrangulations and for large Vorono\"\i\ cells. In particular, both rules above lead to the 
same asymptotic law for the dispatching of area between the two cells.

In this paper, we decide to adopt the first prescription and we define accordingly $F(g,h)$ as the generating function of planar bi-pointed quadrangulation
with a weight $g$ per face strictly within the first Vorono\"\i\ cell, a weight $h$ per face strictly within the second cell and a weight $\sqrt{g\, h}$ per face shared by the two cells. 
Alternatively, $F(g,h)$ is the generating function of i-l.2.f.m with a weight $g$ per 
edge lying strictly in the first face, a weight $h$ per edge lying strictly in the second face and a weight $\sqrt{g\, h}$ per edge incident to both faces.
Our aim will now be to evaluate $F(g,h)$.

\section{Generating function for iso-labelled two-face maps}

\subsection{Connection with the generating function for labelled chains}
\begin{figure}
\begin{center}
\includegraphics[width=9cm]{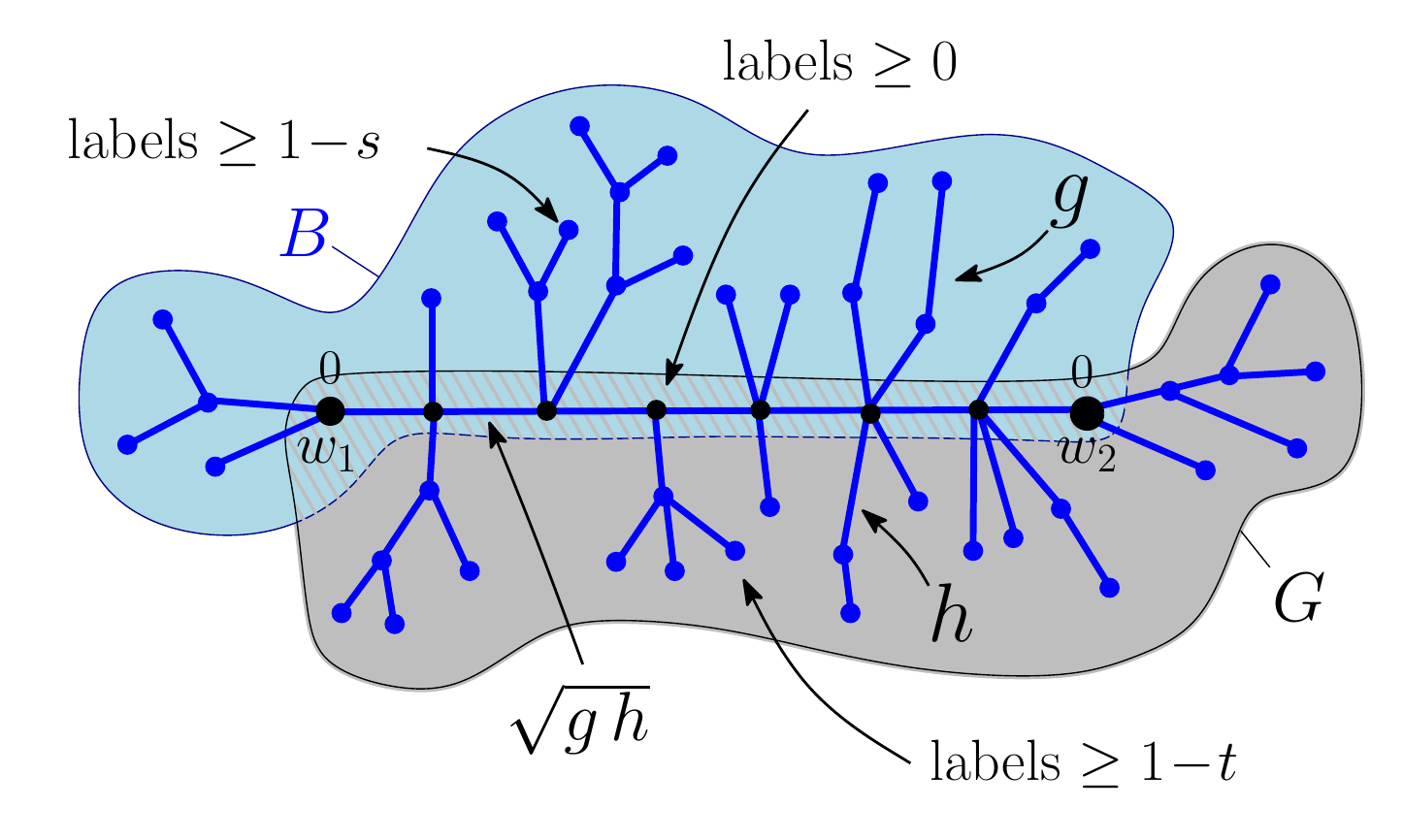}
\end{center}
\caption{A schematic picture of a labelled chain (l.c) contributing to $X_{s,t}(g,h)$ for $s,t>0$. The labels of the vertices in the blue domain $B$ have to be $\geq 1-s$ and
those in the gray domain $G$ have to be $\geq 1-t$. The spine vertices, i.e.\ the vertices of $B\cap G$ are required to have non-negative labels, 
and the spin endpoints $w_1$ and $w_2$ are labelled $0$. The edge weights are: $g$ if they have at least one endpoint in $B\setminus(B\cap G)$, $h$ if they have at least one endpoint in
$G\setminus(B\cap G)$ and $\sqrt{g\, h}$ if they have their two endpoints in  $B\cap G$.}
\label{fig:lc}
\end{figure}
In order to compute $F(g,h)$, let us start by introducing what we call \emph{labelled chains} (l.c), which are planar labelled one-face maps, i.e.\ trees whose vertices carry integer labels
satisfying $\langle \hbox{a}_1\rangle$ and with two distinct (and distinguished) 
marked vertices $w_1$ and $w_2$. Such maps are made of a \emph{spine} which is the unique shortest path in the map joining the two vertices,
naturally oriented from $w_1$ to $w_2$, and of a number of labelled subtrees attached to the spine vertices. All internal (i.e.\ other than $w_1$ and $w_2$) 
spine vertices have two (possibly empty) attached labelled subtrees, one on the left and one on the right. As for $w_1$ and $w_2$, they have a single (possibly empty)
such attached labelled subtree.
For $s$ and $t$ two positive integers, we denote by $X_{s,t}\equiv X_{s,t}(g,h)$ the generating function of planar l.c satisfying (see Figure~\ref{fig:lc}):
\begin{enumerate}[$\langle \hbox{b}_1\rangle$]
\item{$w_1$ and $w_2$ have label $0$. The minimal label for the set of spine vertices is $0$. The edges of the spine receive a weight $\sqrt{g\, h}$.}
\item{The minimal label for the set of vertices belonging to the subtree attached to $w_1$ or to any of the subtrees attached to the left of an internal spine vertex is \emph{larger than or equal to}
 $1-s$.
The edges of these subtrees receive a weight $g$.}
\item{The minimal label for the set of vertices belonging to the subtree attached to $w_2$ or to any of the subtrees attached to the right of an internal spine vertex is  \emph{larger than or equal to} $1-t$.
The edges of these subtrees receive a weight $h$.}
\end{enumerate}
For convenience, we incorporate in $X_{s,t}$ a first additional term $1$ (which may be viewed as the contribution of some ``empty" l.c). For $s,t>0$, we also set $X_{s,0}=X_{0,t}=X_{0,0}=1$.

We now return to the generating function $F(g,h)$ of planar i-l.2.f.m. Let us show that $F(g,h)$ is related to $X_{s,t}$ by the relation:
\begin{equation}
F(g,h)=\sum_{s\geq 1} \Delta_s\Delta_t \log(X_{s,t}(g,h))\Big\vert_{t=s} = \sum_{s\geq 1} \log\left(\frac{X_{s,s}(g,h)X_{s-1,s-1}(g,h)}{X_{s-1,s}(g,h)X_{s,s-1}(g,h)}\right)
\label{eq:FfromX}
\end{equation}
(here $\Delta_s$ denotes the finite difference operator $\Delta_s f(s) \equiv f(s)-f(s-1)$).

As already mentioned, a planar i-l.2.f.m is made of a simple closed loop separating the two faces and with (possibly empty) labelled subtrees attached on both sides of the loop vertices.
The loop may be oriented so as to have the face $f_1$ on its left.
Calling $s$ the minimum label for vertices along the loop, with $s\geq 1$, we may shift all labels by $-s$ and use shifted labels instead of the original ones.
With these shifted labels, the planar i-l.2.f.m enumerated by $F(g,h)$ may alternatively be characterized as follows: there exists positive integers $s$ and $t$ such that:
\begin{enumerate}[$\langle \hbox{c}_1\rangle$]
\item{The minimal label for the set of loop vertices is $0$. The edges of the loop receive a weight $\sqrt{g\, h}$.}
\item{The minimal label for the set of vertices belonging to the subtrees attached to the left of loop vertices (including the loop vertices themselves) is \emph{equal to $1-s$}.
The edges of these subtrees receive a weight $g$.}
\item{The minimal label for the set of vertices belonging to the subtrees attached to the right of loop vertices (including the loop vertices themselves) is \emph{equal to $1-t$}.
The edges of these subtrees receive a weight $h$.}
\item{$s=t$}.
\end{enumerate}
The distinction between $s$ and $t$ might seem somewhat artificial in view of $\langle \hbox{c}_4\rangle$ but it was introduced so that $\langle \hbox{c}_2\rangle$ and $\langle \hbox{c}_3\rangle$
actually mimic the (slightly weaker) constraints $\langle \hbox{b}_2\rangle$ and $\langle \hbox{b}_3\rangle$. 
Returning now to a l.c enumerated by $X_{s,t}(g,h)$, it may, upon cutting the chain at all the internal spine vertices with label $0$, be viewed as a
(possibly empty) \emph{sequence} of an arbitrary number $n\geq 0$ of more restricted l.c whose internal spine vertices all have strictly positive labels, 
enumerated say, by $Z_{s,t}=Z_{s,t}(g,h)$ (with the same edge weights as for $X_{s,t}$). This leads to the simple relation $X_{x,t}=1/ (1-Z_{s,t})$. Similarly, a 
\emph{cyclic sequence} of an arbitrary number $n\geq 1$ of these more restricted l.c is enumerated by $\log(1/ (1-Z_{s,t}))=\log(X_{s,t})$.
For such a cyclic sequence, the concatenation of the spines now forms an oriented loop and $\log(X_{s,t})$ therefore enumerates planar labelled two-face maps with the same characterizations 
as $\langle \hbox{c}_1\rangle$-$\langle \hbox{c}_3\rangle$ above except that the minimum labels on both sides of the loop are now larger than or equal to $1-s$ or $1-t$, 
instead of being exactly equal to $1-s$ and $1-t$. The discrepancy is easily corrected by applying finite difference operators\footnote{Indeed, removing from the set of maps with a minimum label $\geq 1-s$ in $f_1$ those maps with a minimum label $\geq 1-(s-1)=2-s$ amounts
to keeping those maps with minimum label in $f_1$ exactly equal to $1-s$.}, namely by taking instead of  $\log(X_{s,t})$ the function $ \Delta_s\Delta_t \log(X_{s,t})$.
The last requirement $\langle \hbox{c}_4\rangle$ is then trivially enforced by setting $t=s$ in this latter generating function and the summation over the arbitrary value of $s\geq 1$ leads directly to the announced expression \eqref{eq:FfromX}.

The reader will easily check that, as customary in map enumeration problems, the generating function $F(g,h)$ incorporates a symmetry factor $1/k$ for those i-l.2.f.m which display a $k$-fold symmetry\footnote{Maps with
two faces may display a $k$-fold symmetry by rotating them around two ``poles" placed at the centers of the two faces.}. In this paper, we will eventually discuss
results for maps with a large number of edges for which $k$-fold symmetric configurations are negligible.    

\subsection{Recursion relations and known expressions}
Our problem of estimating $F(g,h)$ therefore translates into that of evaluating $X_{s,t}(g,h)$. To this end, we shall need to introduce yet another family of maps,
which are planar one-face labelled maps (i.e\ trees whose vertices carry integer labels
satisfying $\langle \hbox{a}_1\rangle$) which are rooted (i.e.\ have a marked oriented edge), whose root vertex (origin of the root edge) has label $0$ and whose minimal label is  \emph{larger than or equal to} $1-s$,
with $s\geq 1$. We shall denote by $R_s(g)$ ($s\geq 1$) their generating function with a weight $g$ per edge and again, a first term $1$ added for convenience.
This new generating function satisfies the following relation, easily derived by looking at the two subtrees obtained by removing the root edge:
\begin{equation*}
R_s(g)=1+g R_s(g)\left(R_{s-1}(g)+R_s(g)+R_{s+1}(g)\right)
\end{equation*}
for $s\geq 1$, with the convention $R_0(g)=0$. This ``recursion relation" determines $R_s(g)$ for all $s\geq 1$, order by order in $g$.
Its solution was obtained in \cite{GEOD} and reads:
\begin{equation}
R_s(g)=\frac{1+4x+x^2}{1+x+x^2}\frac{(1-x^s)(1-x^{s+3})}{(1-x^{s+1})(1-x^{s+2})}
\quad \hbox{for}\ g= x\frac{1+x+x^2}{(1+4x+x^2)^2}\ .
\label{eq:exactRs}
\end{equation}
Here $x$ is a real in the range $0\leq x\leq 1$, so that $g$ is a real in the range $0\leq g\leq 1/12$. 
Note that the above generating function has a singularity for $g\to 1/12$ even  though the above expression has a well-defined limit for $x\to 1$. 

Knowing  $R_s(g)$, we may easily write down a similar recursion relation for $X_{s,t}(g,h)$, obtained by removing the first edge of the spine: the end
point of this edge either has label $0$ and the remainder of the spine is again a l.c enumerated by $X_{s,t}(g,h)$ or it has label $1$ 
and the remainder of the chain may now be decomposed, by removing the first spine edge leading back to label $0$, into two l.c enumerated by $X_{s+1,t+1}(g,h)$
and  $X_{s,t}(g,h)$ respectively. Extra factors $\sqrt{g\, h}$, $R_s(g)$, $R_t(h)$, $R_{s+1}(g)$ and $R_{t+1}(h)$ are needed to account for the removed edges and 
their attached subtrees (those which are not part of the sub-chains), so that we eventually end up with the relation (see \cite{BG08} for a detailed derivation of this relation
when $h=g$):
\begin{equation}
\hskip -.35cm X_{s,t}(g,h)=1+\sqrt{g\, h}\, R_s(g)R_t(h)X_{s,t}(g,h)\left(\!1\!+\sqrt{g\, h}\, R_{s+1}(g)R_{t+1}(h)X_{s+1,t+1}(g,h)\!\right)
\label{eq:Xstrec}
\end{equation}
valid for non-negative $s$ and $t$. This relation again determines $X_s(g,h)$ for all $s,t\geq 1$ order by order\footnote{By this, we mean that 
$X_{s,t}(\rho g,\rho h)$ is determined order by order in $\rho$.} in $g$ and $h$.

Finding an explicit expression for $X_{s,t}(g,h)$ for arbitrary $g$ and $h$ is a challenging issue which we have not been able to solve. As it will appear, this lack 
of explicit formula is not an unsurmountable obstacle in our quest. Indeed, only the \emph{singularity} of $F(g,h)$ for $g$ and
$h$ tending to their common critical value $1/12$ will eventually matter to enumerate large maps. 
Clearly, the absence of explicit expression or $X_{s,t}(g,h)$ will however make our further discussion much more involved.

Still, we way, as a guideline, rely on the following important result. For $g=h$, an explicit expression for $X_{s,t}(g,h)$ was obtained in \cite{BG08}, namely, for $s,t\geq 0$:
\begin{equation*}
X_{s,t}(g,g)=\frac{(1-x^3)(1-x^{s+1})(1-x^{t+1})(1-x^{s+t+3})}{(1-x)(1-x^{s+3})(1-x^{t+3})(1-x^{s+t+1})}
\quad \hbox{where}\ g= x\frac{1+x+x^2}{(1+4x+x^2)^2}\ .
\end{equation*}

\subsection{Local and scaling limits for the generating functions}
Chapuy's conjecture is for quadrangulations with a \emph{fixed number} $N$ of faces, in the limit of large $N$. Via the Miermont bijection, this corresponds to i-l.2.f.m with
a large fixed number $N$ of edges. Proving the conjecture therefore requires an estimate of the coefficient $[g^{N-\frac{p}{2}}h^\frac{p}{2}]F[g,h]$ (recall that,
due to the weight $\sqrt{g\, h}$ per edge of the loop in i-l.2.f.m, $F(g,h)$ has half integer powers in $g$ and $h$), corresponding to a second Vorono\"\i\ cell of area
$n=p/2$, in the limit of large $N$ and for $\phi=n/N$ of order $1$.
Such estimate in entirely encoded in the singularity of the generating function $F(g,h)$ when the edge weights $g$ and $h$ tend simultaneously to their common singular 
value $1/12$. This leads us to set
\begin{equation}
g=\frac{1}{12}\left(1-\frac{a^4}{36}\epsilon^4\right)\ , \qquad h=\frac{1}{12}\left(1-\frac{b^4}{36}\epsilon^4\right)
\label{eq:scalgh}
\end{equation}
(with a factor $1/36$ and a fourth power in $\epsilon$ for future convenience) and to look at the small $\epsilon$ expansion of $F(g,h)$. 

Before we discuss the case of $F(g,h)$ itself, let us return for a while to the quantities $R_s(g)$ and $X_s(g,g)$ for
which we have explicit expressions. The small $\epsilon$ expansion for $R_s(g)$ may be obtained from \eqref{eq:exactRs} upon
first inverting the relation between $g$ and $x$ so as to get the expansion:
\begin{equation*}
x=1-a\, \epsilon +\frac{a^2 \epsilon ^2}{2}-\frac{5\, a^3 \epsilon ^3}{24}+\frac{a^4 \epsilon ^4}{12}-\frac{13\, a^5 \epsilon ^5}{384}+\frac{a^6 \epsilon ^6}{72}-\frac{157\, a^7 \epsilon ^7}{27648}+\frac{a^8 \epsilon
   ^8}{432}+O(\epsilon ^{9})\ .
\end{equation*}
Inserting this expansion in \eqref{eq:exactRs}, we easily get, for any \emph{finite} $s$:
\begin{equation*}
\begin{split}
R_s(g)& =2-\frac{4}{(s+1) (s+2)}
-\frac{s (s+3) \left(3 s^2+9 s-2\right)\, a^4 \epsilon ^4}{180 (s+1) (s+2)}\\
&\qquad +\frac{s (s+3) \left(5 s^4+30 s^3+59 s^2+42 s+4\right)\, a^6 \epsilon ^6}{7560 (s+1) (s+2)}
+O(\epsilon ^8)\ .\\
\end{split}
\end{equation*}
The most singular term of $R_s(g)$ corresponds to the term of order $\epsilon^6=(216/a^6) (1-12g)^{3/2}$ (the constant term and the term 
proportional to $\epsilon^4= (36/a^4) (1-12g)$ being 
regular) and we immediately deduce the large $N$ estimate:
\begin{equation*}
[g^N]R_s(g) \underset{N \to \infty}{\sim} \frac{3}{4}  \frac{12^{N}}{\sqrt{\pi} N^{5/2}} \frac{s (s+3) \left(5 s^4+30 s^3+59 s^2+42 s+4\right)}{35 (s+1) (s+2)}\ .
\end{equation*}
The above $\epsilon$ expansion corresponds to what is called the \emph{local limit} where $s$ is kept finite when $g$ tends to $1/12$ in $R_s(g)$ (or equivalently when $N\to \infty$
in $[g^N]R_s(g)$).
Another important limit corresponds to the so-called \emph{scaling limit} where we let $s$ tend to infinity when $\epsilon\to 0$ by setting
\begin{equation*}
\hskip 5.cm s=\frac{S}{\epsilon}
\end{equation*}
with $S$ of order $1$.  
Inserting this value in the local limit expansion above, we now get at leading order the expansion 
\begin{equation*}
R_{\left\lfloor S/\epsilon\right\rfloor}(g) =2-\frac{4}{S^2}\epsilon^2
-\frac{a^4\, S^2}{60}\epsilon^2
+\frac{a^6\, S^4}{1512}\epsilon^2
+\cdots 
\end{equation*}
where all the terms but the first term $2$ now contribute to the \emph{the same order} $\epsilon^2$.
This is also the case for \emph{all the higher order terms} of the local limit expansion (which we did not display) and 
a proper re-summation, incorporating all these higher order terms, is thus required. Again, it is easily deduced directly from 
the exact expression \eqref{eq:exactRs} and reads:
\begin{equation}
R_{\left\lfloor S/\epsilon\right\rfloor}(g) =2+r(S,a)\ \epsilon^2+O(\epsilon^3)\ , \qquad r(S,a)=-\frac{a^2 \left(1+10 e^{-a S}+e^{-2 a S}\right)}{3 \left(1-e^{-a S}\right)^2}\ .
 \label{eq:expRg}
 \end{equation}
At this stage, it is interesting to note that the successive terms of the local limit expansion, at leading order in $\epsilon$ for $s=S/\epsilon$, correspond 
precisely to the small $S$ expansion of the scaling function $r(S,a)$, namely:
 \begin{equation*}
r(S,a)=-\frac{4}{S^2}-\frac{a^4\, S^2}{60}+\frac{a^6 S^4}{1512}+O(S^6)\ .
\end{equation*}
In other words, we read from the small $S$ expansion of the scaling function the leading large $s$ behavior of the successive coefficients of
the local limit expansion of the associated generating function.

Similarly, from the exact expression of $X_{s,t}(g,g)$, 
we have the local limit expansion
\begin{equation*}
\begin{split}
\hskip -1.2cm X_{s,t}(g,g) &=3-\frac{6 \left(3+4(s+t)+s^2+ s t +t^2\right)}{(s+3) (t+3) (s+t+1)}-\frac{s (s+1) t (t+1) (s+t+3) (s+t+4) a^4 \, \epsilon ^4}{40 (s+3) (t+3) (s+t+1)}\\
&+ \frac{s (s+1) t (t+1) (s+t+3) (s+t+4) \left(5 (s^2+st+t^2)+20 (s+t)+29\right) a^6 \, \epsilon ^6}{5040 (s+3) (t+3) (s+t+1)}+O(\epsilon ^8)\\
   \end{split}
\end{equation*}
and thus
\begin{equation}
\begin{split}
\Delta_s\Delta_t \log(X_{s,t}(g,g))\Big\vert_{t=s}&=\log \left(\frac{s^2 (2 s+3)}{(s+1)^2 (2 s-1)}\right)-\frac{(2 s+1)a^4 \, \epsilon ^4}{60} 
\\
&\ \ +\frac{(2 s+1) \left(10 s^2+10 s+1\right) a^6 \, \epsilon ^6}{1890}+O(\epsilon ^8)\ .\\
   \end{split}
   \label{eq:logloc}
\end{equation}
Alternatively, we also have the corresponding scaling limit counterparts
\begin{equation}
\begin{split}
& X_{\left\lfloor S/\epsilon\right\rfloor,\left\lfloor T/\epsilon\right\rfloor}(g,g)  =3+x(S,T,a)\ \epsilon+O(\epsilon^2)\ , \\
& \hskip 2.cm x(S,T,a) =-3\, a-\frac{6 a \left(e^{-a S}+e^{-a T}-3 e^{-a (S+T)}+e^{-2a (S+T)}\right)}{\left(1-e^{-a S}\right) \left(1-e^{-a T}\right) \left(1-e^{-a (S+T)}\right)}\\
   \end{split}
   \label{eq:xa}
\end{equation}  
and
\begin{equation}
\begin{split}
\Delta_s\Delta_t \log\left(X_{\left\lfloor S/\epsilon\right\rfloor,\left\lfloor T/\epsilon\right\rfloor}(g,g) \right) \Big\vert_{T=S}&=
\epsilon^2 \partial_S\partial_T \log\left(3+x (S,T,a)\, \epsilon \right)\Big\vert_{T=S}+O(\epsilon^4)\\ &=
\epsilon^3\, \frac{1}{3}\ \partial_S\partial_T x(S,T,a)\Big\vert_{T=S} +O(\epsilon^4)\\
&=\epsilon^3\  \frac{2 \, a^3\, e^{-2 a S} \left(1+e^{-2 a S}\right)}{\left(1-e^{-2 a S}\right)^3} +O(\epsilon^4)\\
&= \epsilon^3\, \left(\frac{1}{2\, S^3}-\frac{a^4 S}{30}+\frac{2 a^6 S^3}{189}+O(S^5)\right)+O(\epsilon^4)\ .\\
\end{split}
\label{eq:logscal}
\end{equation}
Again, we directly read on the small $S$ expansion above the large $s$ leading behaviors of the coefficients
in the local limit expansion \eqref{eq:logloc}. In particular, we have the large $s$ behavior:
\begin{equation*}
\log \left(\frac{s^2 (2 s+3)}{(s+1)^2 (2 s-1)}\right)= \frac{1}{2\, s^3}+O\!\left(\frac{1}{s^4}\right)\ .
\end{equation*}

\subsection{Getting singularities from scaling functions}
We will now discuss how the connection between the local limit and the scaling limit allows us to estimate the dominant singularity of
generating functions of the type of \eqref{eq:FfromX} \emph{from the knowledge of scaling functions only}. As a starter, let us
suppose that we wish to estimate the leading singularity of the 
quantity
\begin{equation}
F(g,g)=\sum_{s\geq 1} \Delta_s\Delta_t \log(X_{s,t}(g,g))\Big\vert_{t=s}
\label{eq:Fgg}
\end{equation}
from the knowledge of $x(S,T,a)$ only.
The existence of the scaling limit allows us to write, for any fixed $S_0$ :
\begin{equation*}
\sum_{s\geq \left\lfloor S_0/\epsilon\right\rfloor} \Delta_s\Delta_t \log(X_{s,t}(g,g))\Big\vert_{t=s}= \epsilon^2 \int_{S_0}^\infty dS\, \frac{1}{3} \partial_S\partial_T x(S,T,a)\Big\vert_{T=S}+O(\epsilon^3)\ .
\end{equation*}
To estimate the missing part in the sum \eqref{eq:Fgg}, corresponding to values of $s$ between 
$1$ and $\left\lfloor S_0/\epsilon\right\rfloor -1$, we recall that the local limit expansion \eqref{eq:logloc}
and its scaling limit counterpart \eqref{eq:logscal}
are intimately related in the sense the we directly read on the small $S$ expansion \eqref{eq:logscal} the large $s$ leading behaviors of the coefficients
in the local limit expansion \eqref{eq:logloc}. More precisely, for $k>0$, the coefficient of $\epsilon^k$ in \eqref{eq:logloc} is a rational function of $s$ 
which behaves at large $s$ like $A_{k-3} s^{k-3}$ where $A_{k-3}$ is the coefficient of $S^{k-3}$
in the small $S$ expansion \eqref{eq:logscal} . Here it is important to note that the allowed values of $k>0$ are even integers starting from $k=4$
(with in particular no $k=2$ term\footnote{If present, this term would give the leading singularity. In its absence, the leading singularity is given by the $\epsilon^6$ term.}). 
Subtracting the $k=0$ term in \eqref{eq:logloc} and \eqref{eq:logscal}, taking the difference and summing over $s$, the above remark implies that
\begin{equation*}
\begin{split}
\hskip -1.2cm \sum_{s=1}^{\left\lfloor S_0/\epsilon\right\rfloor-1}&\left( \left( \Delta_s\Delta_t \log(X_{s,t}(g,g))\Big\vert_{t=s}- \log \left(\frac{s^2 (2 s+3)}{(s+1)^2 (2 s-1)}\right)\right)
\right.\\& - \epsilon^3\, \left.\left(\frac{1}{3} \partial_S\partial_T x(s\, \epsilon,t\, \epsilon,a)\Big\vert_{t=s}-\frac{1}{2\, (s\, \epsilon)^3}\right)\right)= \sum_{s=1}^{\left\lfloor S_0/\epsilon\right\rfloor-1} \sum_{k\geq 4} H_{k-3}(s) \epsilon^k \\
\end{split}
\end{equation*}
where $H_{k-3}(s)$ is a rational function of $s$ which now behaves like $B_{k-3} s^{k-4}$ at large $s$ since the terms of order $s^{k-3}$ cancel out in the difference.
Now, for $k\geq 4$, $\sum_{s=1}^{S_0/\epsilon} H_{k-3}(s)$ behaves for small $\epsilon$ like $B_{k-3}S_0^{k-3}\epsilon^{3-k}/(k-3)$
and the sum above over all terms $k\geq 4$ behaves like $\epsilon^3 \sum_{k\geq 4} B_{k-3} S_0^{k-3}/(k-3)$, hence \emph{is of order $\epsilon^{3}$}.

Since the function $(1/3)\partial_S\partial_T x(S,T,a)\vert_{T=S}-1/(2\, S^3)$ is regular at $S=0$, we may use the approximation 
\begin{equation*}
\hskip -1.2 cm \sum_{s=1}^{\left\lfloor S_0/\epsilon\right\rfloor-1}\!\!\!\! \epsilon^3\, \left(\frac{1}{3} \partial_S\partial_T x(s\, \epsilon,t\, \epsilon,a)\Big\vert_{t=s}\!\!\!\! -\frac{1}{2\, (s\, \epsilon)^3}\right)= \epsilon^2 \int_{\epsilon}^{S_0} \!\! \!\! dS\, \left(\frac{1}{3} \partial_S\partial_T x(S,T,a)\Big\vert_{T=S}\!\!\!\! -\frac{1}{2\, S^3}\right)
+O(\epsilon^3)
\end{equation*}
so that we end up with the estimate
\begin{equation}
\begin{split}
F(g,g)
&= \epsilon^2 \int_{\epsilon}^{\infty} dS\, \frac{1}{3} \partial_S\partial_T x(S,T,a)\Big\vert_{T=S} \\
&+\sum_{s=1}^{\left\lfloor S_0/\epsilon\right\rfloor-1} \log \left(\frac{s^2 (2 s+3)}{(s+1)^2 (2 s-1)}\right)-\epsilon^2 \int_\epsilon^{S_0}dS\, \frac{1}{2S^3}+O(\epsilon^3)\ .\\
\end{split}
\label{eq:expFgg}
\end{equation}
The first term is easily computed to be
\begin{equation*}
\begin{split}
\epsilon^2 \int_{\epsilon}^{\infty} dS\, \frac{1}{3} \partial_S\partial_T x(S,T,a)\Big\vert_{T=S}
& =
\epsilon^2 \int_{\epsilon}^{\infty} dS\, \frac{2 a^3 e^{-2 a S} \left(1+e^{-2 a S}\right)}{\left(1-e^{-2 a S}\right)^3} 
\\& =\epsilon^2  \frac{a^2 e^{-2 a \epsilon }}{\left(1-e^{-2 a \epsilon }\right)^2}= \frac{1}{4}-\frac{a^2 \epsilon ^2}{12}+O(\epsilon ^4)\\
\end{split}
\end{equation*}
and gives us the leading singularity of $F(g,g)$, namely  $-a^2\, \epsilon ^2/12=-(1/2)\sqrt{1-12g}$.
As for the last two terms, their value at small $\epsilon$ is easily evaluated to be
\begin{equation*}
-\frac{1}{4}+\log \left(\frac{4}{3}\right)+O(\epsilon ^4)\ .
\end{equation*}
These terms do not contribute to the leading singularity of $F(g,g)$ and serve only to correct the constant term in the expansion, 
leading eventually to the result:
\begin{equation}
F(g,g)
= \log \left(\frac{4}{3}\right)-\frac{a^2 \epsilon ^2}{12}+O(\epsilon ^3)\ .
\label{eq:singFgg}
\end{equation}

Of course, this result may be verified from the exact expression
\begin{equation*}
\begin{split}
F(g,g)=\sum_{s\geq 1}  \log\left(\frac{X_{s,s}(g,g)X_{s-1,s-1}(g,g)}{X_{s-1,s}(g,g)X_{s,s-1}(g,g)}\right)& =\log\left(\frac{\left(1-x^2\right)^2}{(1-x) \left(1-x^3\right)}\right)
\\ &= \log \left(\frac{4}{3}\right)-\frac{a^2 \epsilon ^2}{12}+O(\epsilon ^3)\\
\end{split}
\end{equation*}
for $x=1-a \epsilon+O(\epsilon^2)$.
The reader might thus find our previous calculation both cumbersome and useless but the lesson of this 
calculation is not the precise result itself but the fact that the leading singularity of a sum like \eqref{eq:Fgg}
is, via \eqref{eq:expFgg}, fully predicable from the knowledge of the scaling function $x(S,T,a)$ only. Note indeed that the singularity is \emph{entirely contained 
in the first term of \eqref{eq:expFgg}} and that the 
last two terms, whose precise form requires the additional knowledge of the first coefficient 
of the local limit of $\Delta_s\Delta_t \log(X_{s,t}(g,g))\vert_{t=s}$ do not contribute to the singularity but serve only to correct the constant term in the expansion which is
not properly captured by the integral of the scaling function. This additional knowledge is therefore not needed strico sensu if we are only interested in 
the singularity of \eqref{eq:Fgg}.

To end this section, we note that we immediately deduce from the leading singularity $-(1/2)\sqrt{1-12g}$ of $F(g,g)$ the large $N$ asymptotics
\begin{equation}
[g^N]F(g,g) \underset{N \to \infty}{\sim} \frac{1}{4}  \frac{12^{N}}{\sqrt{\pi} N^{3/2}}
\label{eq:norm}
\end{equation}
for the number of i-l.2.f.m with $N$ edges, or equivalently, of planar quadrangulations with $N$ faces and with two marked (distinct and distinguished) vertices at even distance from each other.

\section{Scaling functions with two weights $g$ and $h$}
\label{sec:scaling}
\subsection{An expression for the singularity of ${\boldsymbol{F(g,h)}}$}
The above technique gives us a way to access to the singularity of the function $F(g,h)$ via the following small $\epsilon$ estimate, 
which straightforwardly generalizes \eqref{eq:expFgg}:
\begin{equation}
\begin{split}
 F(g,h)
& =\epsilon^2 \int_{\epsilon}^{\infty} dS\,  \frac{1}{3} \partial_S\partial_T x(S,T,a,b)\Big\vert_{T=S} \\
&+\sum_{s=1}^{\left\lfloor S_0/\epsilon\right\rfloor-1} \log \left(\frac{s^2 (2 s+3)}{(s+1)^2 (2 s-1)}\right)-\epsilon^2 \int_\epsilon^{S_0}dS\, \frac{1}{2S^3}+O(\epsilon^3)\ .
\\
\end{split}
\label{eq:expFgh}
\end{equation}
Here $x(S,T,a,b)$ is the scaling function associated to $X_{s,t}(g,h)$ via
\begin{equation}
X_{\left\lfloor S/\epsilon\right\rfloor,\left\lfloor T/\epsilon\right\rfloor}(g,h)  =3+x(S,T,a,b)\ \epsilon+O(\epsilon^2)
\label{eq:xab}
\end{equation}
when $g$ and $h$ tend to $1/12$ as in \eqref{eq:scalgh}. As before, the last two terms of \eqref{eq:expFgh} do not contribute 
to the singularly but give rise only to a constant at this order in the expansion. 
The reader may wonder why these terms are \emph{exactly the same} as those of \eqref{eq:expFgg},
as well as why the leading term $3$ in \eqref{eq:xab} is the same as that of \eqref{eq:xa} 
although $h$ is no longer equal to $g$. This comes from the simple remark that these terms 
all come from the behavior of $X_{s,t}(g,h)$ \emph{exactly at $\epsilon=0$} which is the same as that  
of $X_{s,t}(g,g)$ since, for $\epsilon=0$, both $g$ and $h$ have the same value $1/12$.
In other words, we have 
\begin{equation}
X_{s,t}(g,h)=3-\frac{6 \left(3+4(s+t)+s^2+ s t +t^2\right)}{(s+3) (t+3) (s+t+1)}+O(\epsilon ^4)
 \label{eq:expXgh}
\end{equation}
and consequently, for small $S$ and $T$ of the same order (i.e. $T/S$ finite), we must have an expansion
of the form \eqref{eq:xab} with
\begin{equation}
x(S,T,a,b) =-\frac{6 \left(S^2+ S T +T^2\right)}{S\,T\,(S+T)}+O(S^3)
\label{eq:smalST}
\end{equation}
in order to reproduce the large $s$ and $t$ behavior of the local limit just above. 
We thus have
\begin{equation*}
 \frac{1}{3} \partial_S\partial_T x(S,T,a,b)\Big\vert_{T=S} = \frac{1}{2\, S^3}+O(S)
 \end{equation*}
 while 
 \begin{equation*}
\Delta_s\Delta_t \log(X_{s,t}(g,h))\Big\vert_{t=s}=\log \left(\frac{s^2 (2 s+3)}{(s+1)^2 (2 s-1)}\right)+O(\epsilon ^4)\ ,
\end{equation*}
hence the last two terms in \eqref{eq:expFgh}.
\subsection{An expression for the scaling function ${\boldsymbol {x(S,T,a,b)}}$}
Writing the recursion relation \eqref{eq:Xstrec} for $s=S/\epsilon$ and $t=T/\epsilon$ and using the small $\epsilon$ expansions
\eqref{eq:expRg} and \eqref{eq:xab}, we get at leading order in $\epsilon$ (i.e.\ at order $\epsilon^2$) the following partial differential 
equation\footnote{Here, choosing $(g+h)/2$ instead of $\sqrt{g\, h}$ for the weight of spine edges in the l.c would not
change the differential equation. It can indeed be verified that only the leading value $1/12$ of this weight matters.}
\begin{equation*}
2 \big(x(S,T,a,b)\big)^2+6 \big(\partial_Sx(S,T,a,b)+\partial_Tx(S,T,a,b)\big)+27 \big(r(S,a)+r(T,b)\big)=0
\end{equation*}
which, together with the small $S$ and $T$ behavior \eqref{eq:smalST}, fully determines $x(S,T,a,b)$.
To simplify our formulas, we shall introduce new variables
\begin{equation*}
\sigma\equiv e^{-a S}\ , \qquad \tau\equiv e^{-b T}\ ,
\end{equation*}
together with the associated functions
\begin{equation*}
\mathfrak{X}(\sigma,\tau,a,b)\equiv x(S,T,a,b)\ , \qquad \mathfrak{R}(\sigma,a)\equiv r(S,a)\ .
\end{equation*}
With these variables, the above partial differential equation becomes:
\begin{equation}
\begin{split}
& \hskip -1.2cm 2 \big(\mathfrak{X}(\sigma,\tau,a,b)\big)^2-6 \big(a\, \sigma\, \partial_\sigma\mathfrak{X}(\sigma,\tau,a,b)+b\, \tau\, \partial_\tau\mathfrak{X}(\sigma,\tau,a,b)\big)+27 \big(\mathfrak{R}(\sigma,a)+\mathfrak{R}(\tau,b)\big)=0\\
&\hbox{with}\ \ 
\mathfrak{R}(\sigma,a)=-\frac{a^2 \left(1+10 \sigma +\sigma^2\right)}{3 (1-\sigma)^2}
\ \  \hbox{and}\ \ 
\mathfrak{R}(\tau,b)=-\frac{b^2 \left(1+10 \tau +\tau^2\right)}{3 (1-\tau)^2}\ .
\\
\end{split}
\label{eq:pardif}
\end{equation}
For $b=a$, i.e.\ $h=g$, we already know from \eqref{eq:xa} the solution 
\begin{equation*}
\mathfrak{X}(\sigma,\tau,a,a)=-3a-\frac{6 a \left(\sigma+\tau-3 \sigma \tau + \sigma ^2 \tau ^2 \right)}{(1-\sigma ) (1-\tau) (1-\sigma  \tau)}
\end{equation*}
and it is a simple exercise to check that it satisfies the above partial differential equation in this particular case. This suggests to 
look for a solution of \eqref{eq:pardif} in the form:
\begin{equation*}
\mathfrak{X}(\sigma,\tau,a,b)=-3\sqrt{\frac{a^2+b^2}{2}}-\frac{\mathfrak{N}(\sigma,\tau,a,b)}{(1-\sigma ) (1-\tau) \mathfrak{D}(\sigma,\tau,a,b)}
\end{equation*}
where $\mathfrak{N}(\sigma,\tau,a,b)$ and $\mathfrak{D}(\sigma,\tau,a,b)$ are polynomials in the variables $\sigma$ and $\tau$. The first
constant term is singularized for pure convenience (as it could be incorporated in $\mathfrak{N}$). Its value is chosen by assuming that the function $\mathfrak{X}(\sigma,\tau,a,b)$
is regular for small $\sigma$ and $\tau$ (an assumption which will be indeed verified a posteriori) in which case, from \eqref{eq:pardif},
we expect:
\begin{equation*}
(\mathfrak{X}(0,0,a,b))^2=-\frac{27}{2}\left(\mathfrak{R}(0,a)+\mathfrak{R}(0,b)\right)= 9\, \frac{a^2+b^2}{2}
\end{equation*}
(the $-$ sign is then chosen so as to reproduce the known value $-3a$ for $b=a$).
To test our Ansatz, we tried for $\mathfrak{N}(\sigma,\tau,a,b)$ a polynomial of maximum degree $3$ in $\sigma$ and in $\tau$ 
and for $\mathfrak{D}(\sigma,\tau,a,b)$ a polynomial of maximum degree $2$, namely
\begin{equation*}
\begin{split}
\mathfrak{N}(\sigma,\tau,a,b)&=\sum_{i=0}^{3}\sum_{j=0}^3 n_{i,j}\, \sigma^i\tau^j \ , \\
\mathfrak{D}(\sigma,\tau,a,b)& =\sum_{i=0}^{2}\sum_{j=0}^2 d_{i,j}\, \sigma^i\tau^j \ , \\
\end{split}
\end{equation*}
with $d_{0,0}=1$ (so as to fix the, otherwise arbitrary, normalization of all coefficients, assuming that $d_{0,0}$ does not vanish).
With this particular choice, solving \eqref{eq:pardif} translates, after reducing to the same denominator, into canceling all coefficients of a polynomial
of degree $6$ in $\sigma$ as well as in $\tau$, hence into solving a system of $7\times 7 = 49$ equations for the $4\times 4+3\times 3-1=24$ variables 
$(n_{i,j})_{0\leq i,j\leq 3}$ and $(d_{i,j})_{0\leq i,j\leq 2\atop (i,j)\neq(0,0)}$. Remarkably enough, this system, although clearly
over-determined, \emph{admits a unique solution}
displayed explicitly in Appendix~A. Moreover, we can check from the explicit form of $\mathfrak{N}(\sigma,\tau,a,b)$ and $\mathfrak{D}(\sigma,\tau,a,b)$
the small $S$ and $T$ expansions (with $T/S$ finite):
\begin{equation*}
\begin{split}
& \hskip -1.2cm \mathfrak{N}(e^{-a S},e^{-b T},a,b)= 6\, a\, b\, (S^2+S\, T+T^2)\, \mathfrak{C}(a,b) +O(S^3)\ ,\\
& \hskip -1.2cm \mathfrak{D}(e^{-a S},e^{-b T},a,b)= (S+T)\, \mathfrak{C}(a,b) +O(S^2)\ ,\\
&\hskip -1.2cm \hbox{with}\ \mathfrak{C}(a,b) = \frac{216\, a^2 b^2 \left(a^2+b^2\right) \left(a^2+a\, b+b^2\right)}{(a-b)^2 (a+b) \left(2 a^2+b^2\right) \left(a^2+2 b^2\right)}
 -\frac{36\, \sqrt{2} \, a^2 b^2 \sqrt{a^2+b^2} \left(4\, a^2+a\, b+4\, b^2\right)}{(a-b)^2 \left(2 a^2+b^2\right) \left(a^2+2 b^2\right)} 
\end{split}
\end{equation*}
and (by further pushing the expansion for $\mathfrak{X}$ up to order $S^2$) that 
\begin{equation*}
\mathfrak{X}(e^{-a S},e^{-b T},a,b) =-\frac{6 \left(S^2+ S T +T^2\right)}{S\,T\,(S+T)}+O(S^3)
\end{equation*}
which is the desired initial condition \eqref{eq:smalST}. We thus have at our disposal \emph{an explicit expression} for the
scaling function $\mathfrak{X}(\sigma,\tau,a,b) $, or equivalently $x(S,T,a,b)$ for arbitrary $a$ and $b$.

\subsection{The integration step}
Having an explicit expression for $x(S,T,a,b)$, the next step is to compute the first integral in \eqref{eq:expFgh}.
We have, since setting $T=S$ amounts to setting $\tau=\sigma^{b/a}$:  
\begin{equation*}
\int_{\epsilon}^{\infty} dS\,  \frac{1}{3} \partial_S\partial_T x(S,T,a,b)\Big\vert_{T=S}=
\int_{0}^{e^{-a\, \epsilon}} d\sigma\,  \frac{1}{3}\, b\, \sigma^{b/a} \partial_\sigma\partial_\tau \mathfrak{X}(\sigma,\tau,a,b)\Big\vert_{\tau=\sigma^{b/a}}\ .
\end{equation*}
To compute this latter integral, it is sufficient to find a primitive of its integrand, namely a function $\mathfrak{K}(\sigma,a,b)$ such that:
\begin{equation}
\partial_\sigma \mathfrak{K}(\sigma,a,b) = \frac{1}{3}\, b\, \sigma^{b/a} \partial_\sigma\partial_\tau \mathfrak{X}(\sigma,\tau,a,b)\Big\vert_{\tau=\sigma^{b/a}}\ .
\label{eq:K}
\end{equation}
For $b=a$, we have from the explicit expression of $\mathfrak{X}(\sigma,\tau,a,a)$:
\begin{equation*}
\begin{split}
\hskip -1.cm  \frac{1}{3}\, a\, \sigma \partial_\sigma\partial_\tau \mathfrak{X}(\sigma,\tau,a,a)\Big\vert_{\tau=\sigma}= \frac{2 a^2 \sigma  \left(1+\sigma ^2\right)}{\left(1-\sigma ^2\right)^3} & =\partial_\sigma 
 \left(\frac{a^2 \sigma ^2}{\left(1-\sigma ^2\right)^2}\right)\\
 & 
 = \partial_\sigma \left(\frac{a\, b \, \sigma\, \tau }{\left(1-\sigma\, \tau\right)^2}\Bigg\vert_{\tau=\sigma^{b/a}}\right)\Bigg\vert_{b=a}\ .\\
 \end{split}
\end{equation*}
In the last expression, we recognize the square of the last factor $(1-\sigma\, \tau)$ appearing in the denominator in $\mathfrak{X}(\sigma,\tau,a,a)$. 
This factor is replaced by $\mathfrak{D}(\sigma,\tau,a,b)$ when $b\neq a$ and this suggest to look for 
an expression of the form:
\begin{equation*}
\mathfrak{K}(\sigma,a,b)=\frac{a\, b\, \sigma\, \tau\, \mathfrak{H}(\sigma,\tau,a,b)}{\left(\mathfrak{D}(\sigma,\tau,a,b)\right)^2}\Bigg\vert_{\tau=\sigma^{b/a}}
\end{equation*}
with the same function $\mathfrak{D}(\sigma,\tau,a,b)$ as before and where $\mathfrak{H}(\sigma,\tau,a,b)$ is now a polynomial of the form
\begin{equation*}
\mathfrak{H}(\sigma,\tau,a,b)=\sum_{i=0}^{2}\sum_{j=0}^2 h_{i,j}\, \sigma^i\tau^j
\end{equation*}
(here again the degree $2$ in each variable $\sigma$ and $\tau$ is a pure guess).
With this Ansatz, eq.~\eqref{eq:K} translates, after some elementary manipulations, into
\begin{equation*}
 \frac{1}{3}\partial_\sigma\partial_\tau \mathfrak{X}(\sigma,\tau,a,b)= \left\{(a+b)+a\, \sigma\, \partial_\sigma +b\, \tau\, \partial_\tau \right\}\frac{\mathfrak{H}(\sigma,\tau,a,b)}{\left(\mathfrak{D}(\sigma,\tau,a,b)\right)^2}
 \end{equation*}
 which needs being satisfied only for $\tau=\sigma^{b/a}$. We may however decide to look for a function $\mathfrak{H}(\sigma,\tau,a,b)$ which 
 satisfies the above requirement for arbitrary independent values of $\sigma$ and $\tau$. After reducing to the same denominator, 
 we again have to cancel the coefficients of a polynomial of degree $6$ in $\sigma$ as well as in $\tau$. This gives rise to a system of $7\times 7=49$ equations for the 
 $3\times 3=9$ variables $(h_{i,j})_{0\leq i,j\leq 2}$. Remarkably enough, this over-determined system
 again admits a unique solution displayed explicitly in Appendix~B.
 
This solution has non-zero finite values for $\mathfrak{H}(0,0,a,b)$ and $\mathfrak{D}(0,0,a,b)$ and therefore we deduce $\mathfrak{K}(0,a,b)=0$ so that we find
\begin{equation*}
\begin{split}
\hskip -1.2cm \int_{0}^{e^{-a\, \epsilon}} d\sigma\,  \frac{1}{3}\, b\, \sigma^{b/a} \partial_\sigma\partial_\tau \mathfrak{X}(\sigma,\tau,a,b)\Big\vert_{\tau=\sigma^{b/a}}
&=\mathfrak{K}(e^{-a\, \epsilon},a,b)\\
&= \frac{a\, b\, e^{-a\, \epsilon}\, e^{-b\, \epsilon}\, \mathfrak{H}(e^{-a\, \epsilon},e^{-b\, \epsilon},a,b)}{\left(\mathfrak{D}(e^{-a\, \epsilon},e^{-b\, \epsilon},a,b)\right)^2}\\
&=\frac{1}{4\, \epsilon^2}-\frac{\left(a^2- a\, b+b^2\right) \left(a^2+a\, b+b^2\right)}{18 \left(a^2+b^2\right)}+O(\epsilon^2)\ .
\end{split}
\end{equation*}
Eq.~\eqref{eq:expFgh} gives us the desired singularity
\begin{equation}
\begin{split}
 F(g,h)
&=\frac{1}{4}-\frac{\left(a^2-a\, b+b^2\right) \left(a^2+a\, b+b^2\right)}{18 \left(a^2+b^2\right)}\epsilon^2
-\frac{1}{4}+\log \left(\frac{4}{3}\right)+O(\epsilon ^3)\\
&= \log \left(\frac{4}{3}\right) -\frac{\left(a^2-a\, b+b^2\right) \left(a^2+a\, b+b^2\right)}{18 \left(a^2+b^2\right)}\epsilon^2+O(\epsilon ^3)\\
&= \log \left(\frac{4}{3}\right) -\frac{1}{18}\, \frac{(a^6-b^6)}{(a^4-b^4)}\, \epsilon^2+O(\epsilon ^3)\ .\\
\end{split}
\label{eq:lutfin}
\end{equation}
Note that for $b=a$ ($h=g$), we recover the result \eqref{eq:singFgg} for the singularity of $F(g,g)$, as it should. 

More interestingly, 
we may now obtain from \eqref{eq:lutfin} some asymptotic estimate for the number $[g^{N-\frac{p}{2}}h^{\frac{p}{2}}]F(g,h)$ of planar quadrangulations with $N$ faces, with two marked (distinct and distinguished) vertices at even distance from each other and with Vorono\"\i\ cells of respective areas $N-(p/2)$ and $(p/2)$ (recall
that, due to the existence of faces shared by the two cells, the area of a cell may be any half-integer between $0$ and $N$).
Writing
\begin{equation}
\begin{split}
\hskip -1.2cm -\frac{1}{18}\, \frac{(a^6-b^6)}{(a^4-b^4)}\, \epsilon^2&=-\frac{1}{18}\,  \frac{(a^4\, \epsilon^4)^{3/2}-(b^4\, \epsilon^4)^{3/2}}{(a^4\, \epsilon^4)-(b^4\, \epsilon^4)}
\\&=\frac{1}{36}\, \frac{(1-12h)^{3/2} -(1-12g)^{3/2}}{h-g}
\\&=\frac{1}{6}\, \frac{\sqrt{h}(1-12h)^{3/2} -\sqrt{g}(1-12g)^{3/2}}{\sqrt{h}-\sqrt{g}}+O(\epsilon^6)
\\&=\frac{1}{6}+\sum_{N\geq 1} \frac{2N\!+\!1}{2N\!-\!3}\, \frac{3^N}{N}{2(N\!-\!1)\choose N\!-\!1}\sum_{p=0}^{2N} \frac{g^{N-\frac{p}{2}}\, h^{\frac{p}{2}}}{2N+1} +O(\epsilon^6)\ ,
\end{split}
\label{eq:ident}
\end{equation}
where we have on purpose chosen in the third line an expression whose expansion involves half integer powers in $g$ and $h$, we deduce \emph{heuristically} that, \emph{for large $N$}, $[g^{N-\frac{p}{2}}h^{\frac{p}{2}}]F(g,h)$ behaves like
\begin{equation*}
[g^{N-\frac{p}{2}}h^{\frac{p}{2}}] \frac{1}{6}\, \frac{\sqrt{h}(1-12h)^{3/2} -\sqrt{g}(1-12g)^{3/2}}{\sqrt{h}-\sqrt{g}}
\underset{N \to \infty}{\sim} \frac{1}{4}  \frac{12^{N}}{\sqrt{\pi} N^{3/2}} \times  \frac{1}{2N+1}
\end{equation*}
\emph{independently of $p$}. After normalizing by \eqref{eq:norm}, the probability that the second Vorono\"\i\ cell has some fixed half-integer area $n=p/2$ ($0\leq p\leq 2N$) is asymptotically 
equal to $1/(2N+1)$ independently
of the value of $n$. As a consequence, \emph{the law for $\phi=n/N$ is uniform in the interval $[0,1]$.}

Clearly, the above estimate is too precise and has no reason to be true {\it{stricto sensu}} for finite values of $p$. Indeed, in the expansion \eqref{eq:lutfin}, both $g$ and $h$ tend simultaneously to
$0$, so that the above estimate for $[g^{N-\frac{p}{2}}h^{\frac{p}{2}}]F[g,h]$ should be considered as valid only when both $N$ and $n=p/2$ become large in a limit where the ratio $\phi=n/N$ may
be considered as a finite continuous variable. In other word, some \emph{average} over values of $p$ with $n=p/2$ in the range $N\phi \leq n< N(\phi+d\phi)$ is implicitly required. 
With this averaging procedure, any other generating function with the same singularity as \eqref{eq:lutfin}  would then lead to the same uniform law for $\phi$.
For instance, using the second line of \eqref{eq:ident} and writing
\begin{equation*}
\frac{1}{36}\, \frac{(1-12h)^{3/2} -(1-12g)^{3/2}}{h-g}
= -\frac{1}{2}+\sum_{N\geq 1} \frac{3^N}{N}{2(N\!-\!1)\choose N\!-\!1}\sum_{n=0}^N \frac{g^{N-n}\, h^n}{N+1}\ ,
\end{equation*}
we could as well have estimated from our singularity a value of $[g^{N-\frac{p}{2}}h^{\frac{p}{2}}]F(g,h)$ asymptotically equal to:
\begin{equation*}
[g^{N-\frac{p}{2}}h^{\frac{p}{2}}]\frac{1}{36}\, \frac{(1-12h)^{3/2} -(1-12g)^{3/2}}{h-g}
\underset{N \to \infty}{\sim} \frac{1}{4}  \frac{12^{N}}{\sqrt{\pi} N^{3/2}} \times  \frac{1}{N+1} \, \delta_{p,\rm{even}}
\end{equation*}
with $\delta_{p,\rm{even}}=1$ if $p$ is even and $0$ otherwise. Of course, averaging over both parities, this latter estimate leads to the same uniform 
law for the continuous variable $\phi=n/N$.
\vskip .2cm
Beyond the above heuristic argument, we may compute the law for $\phi$ in a rigorous way by considering the  large $N$ behavior of the fixed $N$ expectation value 
\begin{equation*}
E_N[e^{\mu\, \left(\frac{n}{N}\right)}]\equiv\frac{\sum\limits_{p=0}^{2N} e^{\mu\, \left(\frac{p}{2N}\right)}\, [g^{N-\frac{p}{2}} h^{\frac{p}{2}}]F[g,h]}{\sum\limits_{p=0}^{2N} [g^{N-\frac{p}{2}} 
h^{\frac{p}{2}}]F[g,h]}
=\frac{[g^N]F(g,g\, e^{\frac{\mu}{N}})}{[g^N]F(g,g)}\ .
\end{equation*}
The coefficient $[g^N]F(g,g\, e^{\frac{\mu}{N}})$ may then be obtained by a contour integral around $g=0$, namely
\begin{equation*}
\frac{1}{2\rm{i}\pi} \oint\frac{dg}{g^{N+1}} F(g,g\, e^{\frac{\mu}{N}})
\end{equation*}
and, at large $N$, we may use \eqref{eq:scalgh} and \eqref{eq:lutfin} with 
\begin{equation*}
\hskip 5.cm \epsilon^4=\frac{1}{N}
\end{equation*}
 to rewrite this integral as an integral over $a$. More precisely, at leading order in $N$, setting $h=g\, e^{\frac{\mu}{N}} =g\, (1+ \mu \epsilon^4)$
amounts to take:
\begin{equation*}
\hskip 4.5cm b^4=a^4-36\mu\ . 
\end{equation*}
Using $dg=-(1/12)a^3/(9 N)$, $g^{N+1}\sim (1/12)^{N+1}\, e^{-a^4/36}$ (and ignoring the constant term $\log(4/3)$ which does not contribute to 
the $g^N$ coefficient for $N\geq 1$), the contour integral above becomes at leading order:
\begin{equation*}
\frac{1}{2\rm{i}\pi} \frac{12^N}{N^{3/2}} \int_\mathcal{C} da\,  \frac{-a^3}{9} \left\{-\frac{1}{18} \frac{a^6-(a^4-36\mu)^{3/2})}{36\mu}\, e^{a^4/36} +O\!\left(\frac{1}{N^{1/4}}\right)\right\}
\end{equation*}
where the integration path follows some appropriate contour $\mathcal{C}$ in the complex plane. The precise form of this contour and the details of the computation of this integral are given in Appendix~C. 
We find the value\begin{equation*}
\frac{1}{2\rm{i}\pi} \int_\mathcal{C} da\,  \frac{-a^3}{9} \left\{ -\frac{1}{18} \frac{a^6-(a^4-36\mu)^{3/2})}{36\mu}\, e^{a^4/36} \right\}=\frac{1}{4 \sqrt{\pi}} \times  \frac{e^\mu-1}{\mu}\ ,
\end{equation*}
which matches the asymptotic result obtained by
the identification \eqref{eq:ident} since\footnote{We have as well:
$$\frac{1}{4}  \frac{12^{N}}{\sqrt{\pi} N^{3/2}} \times \sum_{p=0}^{2N} \frac{1}{N+1} e^{\mu\, \left(\frac{p}{2N}\right)}\, \delta_{p,\rm{even}}=\frac{1}{4}  \frac{12^{N}}{\sqrt{\pi} N^{3/2}} \times
\sum_{n=0}^{N}\frac{1}{N+1}e^{\mu\, \left(\frac{n}{N}\right)}
 \underset{N \to \infty}{\sim} \frac{1}{4}  \frac{12^{N}}{\sqrt{\pi} N^{3/2}} \times  \frac{e^\mu-1}{\mu}\ .$$}
\begin{equation*}
\frac{1}{4}  \frac{12^{N}}{\sqrt{\pi} N^{3/2}} \times \sum_{p=0}^{2N} \frac{1}{2\,N+1} e^{\mu\, \left(\frac{p}{2N}\right)} \underset{N \to \infty}{\sim} \frac{1}{4}  \frac{12^{N}}{\sqrt{\pi} N^{3/2}} \times  \frac{e^\mu-1}{\mu}\ .
\end{equation*}
After normalization by $[g^N]F(g,g)$ via \eqref{eq:norm}, we end up with the result
\begin{equation*}
\hskip 3.cm E_N[e^{\mu\, \left(\frac{n}{N}\right)}]\underset{N \to \infty}{\sim}\frac{e^\mu-1}{\mu}\ .
\end{equation*}
Writing
\begin{equation*}
\hskip 3.cm \frac{e^\mu-1}{\mu}\ = \int_0^1 d\phi\, e^{\mu\, \phi} \, \mathcal{P}(\phi) \ ,
\end{equation*}
where $\mathcal{P}(\phi)$ is the law for the proportion of area $\phi=n/N$ in, say the second Vorono\"\i\ cell, we obtain that
\begin{equation*}
\hskip 3.cm \mathcal{P}(\phi)=1\quad \forall \phi\in[0,1]\ ,
\end{equation*} 
i.e.\ the law is uniform on the unit segment. This proves the desired result and corroborates Chapuy's conjecture. 
\vskip .2cm
To end our discussion on quadrangulations, let us mention a way to extend our analysis to the case where the distance $d(v_1,v_2)$ is
equal to some odd integer. Assuming that this integer is at least $3$, we can still use the Miermont bijection at the price of introducing a ``delay" $1$ 
for one of two vertices, namely labelling now the vertices by, for instance, $\ell(v)=\min(d(v,v_1), d(v, v_2) + 1)$ and repeating the construction of Figure~\ref{fig:SR}. 
This leads to a second Vorono\"\i\ cell slightly smaller (on average) than the first one but this effect can easily be corrected by averaging the law for $n$ and that for $N-n$. At large 
$N$, it is easily verified that the generating function generalizing $F(g,h)$ to this (symmetrized) ``odd" case (i.e.\ summing over all values $d(v_1,v_2)=2s+1$, $s\geq 1$) has a similar expansion as \eqref{eq:lutfin}, except for the constant term 
$\log(4/3)$ which is replaced by the different value $\log(9/8)$. What matters however is that this new generating function has the same singularity as before when $g$ and 
$h$ tend to $1/12$ so that we still get the uniform law $\mathcal{P}(\phi)$ for the ratio $\phi= n/N$ at large $N$. Clearly, summing over both parities of $d(v_1,v_2)$ would then also 
lead to the uniform law for $\phi$.

\section{Vorono\"\i\ cells for general maps}
\subsection{Coding of general bi-pointed maps by i-l.2.f.m}
\begin{figure}
\begin{center}
\includegraphics[width=8cm]{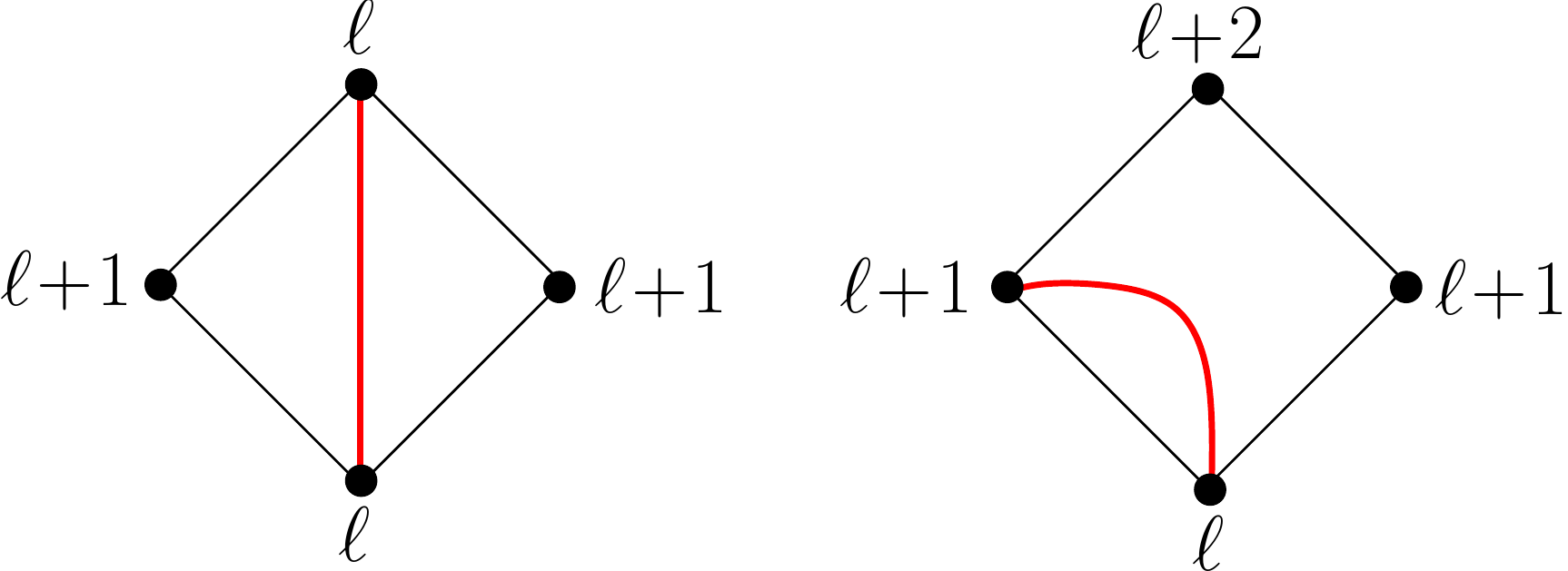}
\end{center}
\caption{The local rules of the Ambj\o rn-Budd bijection.}
\label{fig:AB}
\end{figure}
\begin{figure}
\begin{center}
\includegraphics[width=8cm]{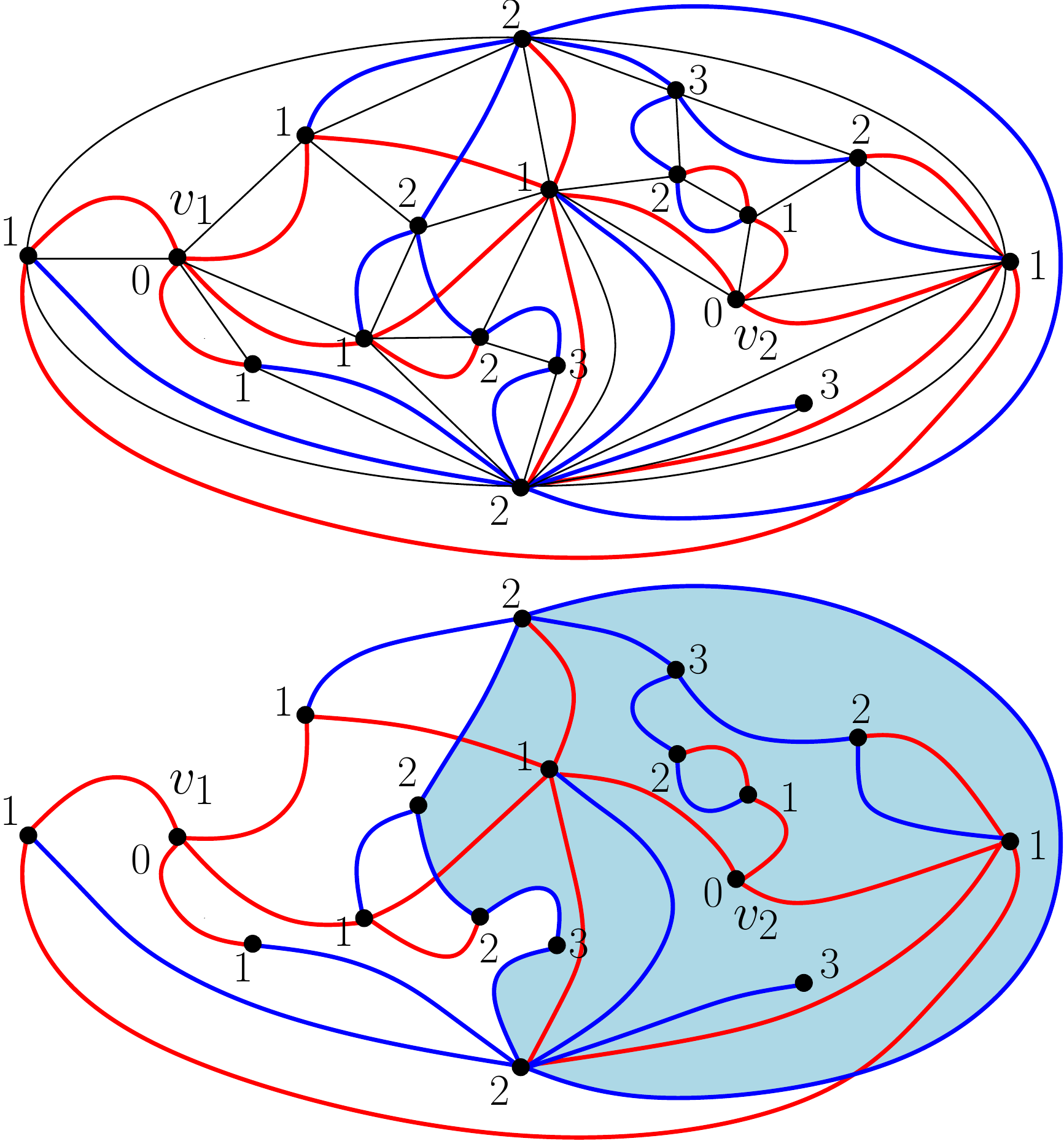}
\end{center}
\caption{Lower part: a general bi-pointed planar map (in red) and the associated i-l.2.f.m (in blue). Upper part: Both maps result from the same bi-pointed quadrangulation using,
on one hand the Miermont bijection via the rules of Figure~\ref{fig:SR} and on the other hand the Ambj\o rn-Budd bijection via the rules of Figure~\ref{fig:AB}.
Note that the label of a vertex $v$ of the general map corresponds to $\min(\delta(v,v_1),\delta(v,v_2))$ where $\delta$ is the graph distance in this map.}
\label{fig:genmap}
\end{figure}
Another direct application of our calculation concerns the statistics of Vorono\"\i\ cells in \emph{bi-pointed general planar maps}. i.e.\ maps with faces of arbitrary degree
and with two distinct (and distinguished) vertices $v_1$ and $v_2$, now at \emph{arbitrary distance $\delta(v_1,v_2)\geq 1$}. As customary, the ``area" of general maps is measured 
by their number $N$ of \emph{edges} to ensure the existence of a finite number of maps for a fixed $N$. General maps are known to be bijectively related
to quadrangulations and it is therefore not surprising that bi-pointed general planar maps may also be coded by i-l.2.f.m. Such a coding is displayed in Figure~\ref{fig:genmap}
and its implementation was first discussed in \cite{AmBudd}.
The simplest way to understand it is to start from a bi-pointed quadrangulation like that of Figure~\ref{fig:cells} (with its two marked vertices $v_1$ and $v_2$
and the induced labelling $\ell(v)=\min(d(v,v_1),d(v,v_2)$) and to draw within each face a new edge according to the rules of 
Figure~\ref{fig:AB} which may be viewed as complementary to the rules of Figure~\ref{fig:SR}. The resulting map formed by these new edges 
is now a general planar map (with faces of arbitrary degree) which is still bi-pointed since 
$v_1$ and $v_2$ are now retained in this map, with vertices labelled by $\ell(v)=\min(\delta(v,v_1),\delta(v,v_2))$ where $\delta(v,v')$ is the graph distance between 
$v$ and $v'$ \emph{in the resulting map}\footnote{Note that, although related, the distance $\delta(v,v')$ between two vertices $v$ and $v'$ in the resulting map and that, $d(v,v')$, in the original quadrangulation are not identical in general.}. This result was shown by Ambj\o rn and Budd in \cite{AmBudd} who also proved that this new construction provides a bijection 
between bi-pointed planar maps with $N$ edges and their two marked vertices at \emph{arbitrary graph distance} and bi-pointed planar quadrangulations with $N$ faces and 
their two marked vertices at \emph{even graph distance}\footnote{In their paper, Ambj\o rn and Budd considered quadrangulations with general labellings satisfying
$\ell(v)-\ell(v')=\pm 1$ if $v$ and $v'$ are adjacent. The present bijection is a specialization of their bijection when the labelling has exactly two local minima (the marked vertices) and the label
is $0$ for both minima. This implies that the two minima are at even distance from each other in the quadrangulation.}.  Note that, in the bi-pointed general map, the labelling may be
erased without loss of information since it may we retrieved directly from graph distances. 
Combined with the Miermont bijection, the Ambj\o rn-Budd bijection gives the desired coding of bi-pointed general planar maps by i-l.2.f.m, whose two faces $f_1$ and $f_2$ moreover surround the
vertices $v_1$ and $v_2$ respectively.
In this coding, all the vertices of the general maps except $v_1$ and $v_2$
are recovered in the i-l.2.f.m, with the same label but the i-l.2.f.m has a number of additional vertices, one lying in each face of the general map and carrying a label 
equal to $1$ plus the maximal label in this face. As discussed in \cite{FG14b}, if the distance $\delta(v_1,v_2)$ is even, equal to $2s $ ($s\geq 1$), 
the i-l.2.f.m (which has by definition minimal label $1$ in its two faces) has a minimum label equal to $s$ for the vertices along the loop separating the two faces, 
and \emph{none of the loop edges has labels $s\, \rule[1.5pt]{8.pt}{1.pt}\, s$}. If the distance $\delta(v_1,v_2)$ is odd, equal to $2s-1$ ($s\geq 1$), 
the i-l.2.f.m has again a minimum label equal to $s$ for the vertices along the loop separating the two faces, but now has  \emph{at least one loop edge with labels $s\, \rule[1.5pt]{8.pt}{1.pt}\, s$}. 

\subsection{Definition of Vorono\"\i\ cells for general maps}
\begin{figure}
\begin{center}
\includegraphics[width=8cm]{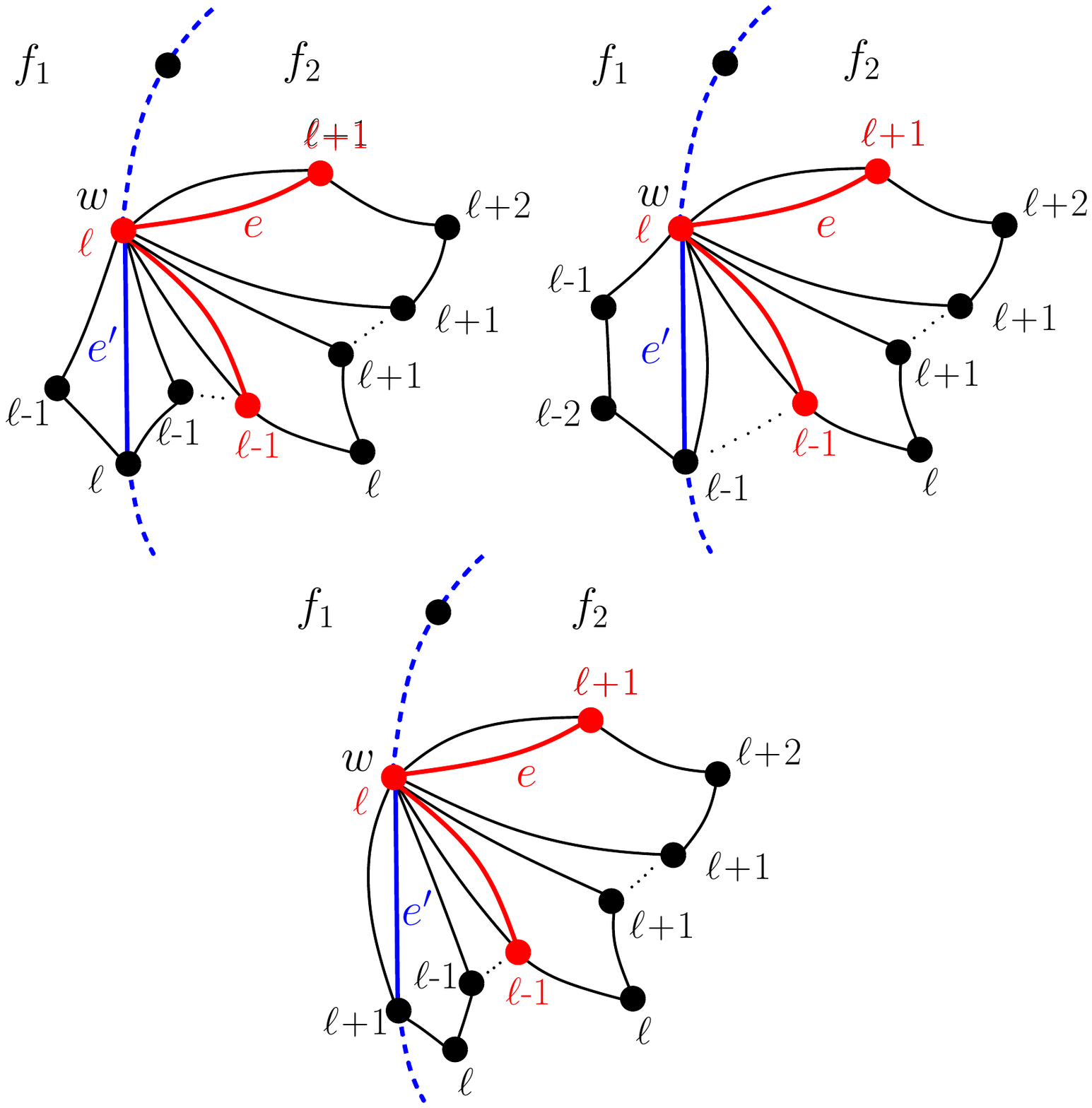}
\end{center}
\caption{Explanation of the ``rebound" property: for any edge $e$ of type $\ell+1\to \ell$ of the general map (here in red) lying in $f_2$ and hitting a loop vertex $w$ of the associated 
i-l.2.f.m (here in blue), there exists, in the sector around $w$ going clockwise from $e$ to the loop edge $e'$ of the i-l.2.f.m leading to $w$ and having $f_1$ on its left, an edge of the general map leaving $w$ within $f_2$ and with endpoint of label $\ell-1$. To see that, we first note that, in the associated quadrangulation (in black),  the first edge leaving $w$ clockwise in the 
sector has label $\ell+1$ from the rules of Figure~\ref{fig:AB}. Similarly the last edge in the quadrangulation leaving $w$ clockwise in the 
sector has label $\ell-1$ from the rules of Figure~\ref{fig:SR}. This holds for the three possible values of the label at the origin of $e'$, namely $\ell$ (upper left), $\ell-1$ (upper right)
and $\ell+1$ (bottom). Then there must be around $w$ in this sector two clockwise consecutive edges of the quadrangulation with
respective labels $\ell+1$ and $\ell-1$ at their endpoint other than $w$. From the rules of Figure~\ref{fig:AB}, the incident face in the quadrangulation gives rise to an edge of the general map lying in the sector, hence in $f_2$, and leaving $w$ toward a vertex with label $\ell-1$. This ``rebound" property is easily generalized to the case where the hitting edge $e$ is of type $\ell\to\ell$. In that case, we only need that the
second half of $e$ lies in $f_2$ to ensure the existence of a subsequent edge $\ell\to\ell-1$ in $f_2$.
}
\label{fig:bump}
\end{figure}

As before, we may define the two Vorono\"\i\ cells in bi-pointed general planar maps as the domains obtained by cutting them along the loop of the associated i-l.2.f.m. 
Let us now see why this definition again matches what we expect from a Vorono\"\i\ cell, namely that vertices in one cell are closer to one of the marked vertices than
to the other. Let us show that any vertex $v$ of the general map strictly inside, say the second face $f_2$ (that containing $v_2$) is closer to $v_2$ than to $v_1$ (or possibly at the same distance).
Since this is obviously true for $v_2$, we may assume $v\neq v_2$ in which case $\ell(v)>0$. Recall that, for any $v$, $\ell(v)=\min(\delta(v,v_1),\delta(v,v_2))$ so that the vertex $v$ necessarily has  a neighbor with label $\ell(v)-1$ within the general map, which itself, if $\ell(v)>1$, has a neighbor of label $\ell(v)-2$, and so on. 
A sequence of edges connecting these neighboring vertices
with strictly decreasing labels provides a shortest path from $v$ to a vertex with label $0$, i.e.\ to either $v_1$ or $v_2$. Let us show
this path may always be chosen so as to stay inside $f_2$, so that it necessarily ends at
$v_2$ and thus $\ell(v)=\delta(v,v_2)\leq \delta(v,v_1)$. 
To prove this, we first note that, since by construction the map edges (in red in the figures) and the i-l.2.f.m edges (in blue) cross only along red edges of type $m\, \rule[1.5pt]{8.pt}{1.pt}\, m$ which cannot belong to a path with strictly decreasing labels, if such a path (which starts with an edge in $f_2$) crosses the loop a first time so as to enter $f_1$, it has to first hit the loop separating the two faces 
at some loop vertex $w$ with, say label $\ell$. We may then rely on the following ``rebound" property, explained in Figure~\ref{fig:bump}:
looking at the environment of $w$ in the sector going clockwise from the map edge $(\ell+1)\to \ell$ of the strictly decreasing path leading
to $w$ (this edge lies in $f_2$ by definition) and the loop edge of the i-l.2.f.m leading to $w$ (with the loop oriented as before with $f_1$ on its left), we see that there always exist a map 
edge $\ell\to \ell-1$ leaving $w$ and lying inside this sector and therefore in $f_2$ (see the legend of Figure~\ref{fig:bump} for
a more detailed explanation). We may then decide to take this edge as the next edge in our path with decreasing labels which de facto, may always be chosen to as to stay\footnote{Note that
some of the vertices along the path may lie on the loop but the path must eventually enter strictly inside $f_2$ since loop labels are larger than $1$.} in $f_2$. 

Let us now discuss vertices of the general map which belong to both Vorono\"\i\ cells, i.e.\ are loop vertices in the i-l.2.f.m. Such vertices may be strictly closer to $v_1$, strictly closer to $v_2$ or
at equal distance from both. More precisely, if a loop vertex $v$ with label $\ell$ is incident to a general map edge in $f_2$, then we can find a path with decreasing labels 
staying inside $f_2$ and thus $\delta(v,v_2)\leq \delta(v,v_1)$. Indeed, if the incident edge is of type $\ell\, \rule[1.5pt]{8.pt}{1.pt}\,(\ell-1)$, it gives
the first step of the desired path, if it is of type $\ell\, \rule[1.5pt]{8.pt}{1.pt}\,(\ell+1)$, looking at this edge backwards and using the rebound property,
the loop vertex is also incident to an edge of type $\ell\, \rule[1.5pt]{8.pt}{1.pt}\,(\ell-1)$ in $f_2$ which may serve as the first step of the desired path. If the incident edge is 
of type $\ell\, \rule[1.5pt]{8.pt}{1.pt}\,\ell$, a straightforward extension of the rebound property shows that the loop vertex is again incident to an edge of 
type $\ell\, \rule[1.5pt]{8.pt}{1.pt}\,(\ell-1)$ in $f_2$ which provides the first step of the desired path.

Similarly, if a loop vertex $v$  incident to a general map edge in $f_1$, then $\delta(v,v_1)\leq \delta(v,v_2)$ and, 
as a consequence, if a loop vertex $v$  in incident to a general map edge in both $f_1$ and $f_2$, then $\delta(v,v_1)= \delta(v,v_2)$.  
From the above properties, we immediately deduce that all the map edges inside $f_1$ (respectively $f_2$) have their two endpoints closer to $v_1$ than to $v_2$ 
(respectively closer to $v_2$ than to $v_1$) or possibly at the same distance.  
As for map edges shared by the two cells, they necessarily connect two vertices $w_1$ and $w_2$ (lying in $f_1$ and $f_2$ respectively) with the same label and with $w_1$ closer to $v_1$ than to $v_2$
(or at the same distance) and $w_2$ closer to $v_2$ than to $v_1$ (or at the same distance).
This fully justifies our definition of Vorono\"\i\ cells.

\subsection{Generating functions and uniform law}
In the context of general maps, a proper measure of the ``area" of Vorono\"\i\ cells is now provided by the number of edges of the general map lying within each cell. Again, 
a number of these edges are actually shared by the two cells, hence contribute $1/2$ to the area of each cell. In terms of generating functions, 
edges inside the first cell receive accordingly the weight $g$, those in the second cell the weight $h$ and those shared by the two cells the weight $\sqrt{g\, h}$
and we  call $F^{\rm{even}}(g,h)$ and $F^{\rm{odd}}(g,h)$ the corresponding generating functions for bi-pointed maps conditioned to have their marked
vertices at even and odd distance respectively.  

When transposed to the associated i-l.2.f.m, this amounts as before to assigning the weight $g$ to those edges of the i-l.2.f.m strictly in $f_1$, $h$ to those strictly in $f_2$, and $\sqrt{g\, h}$ to 
those on the loop separating $f_1$ and $f_2$. Indeed, from the rules of figures \ref{fig:SR} and \ref{fig:AB}, edges of the i-l.2.f.m are in one-to-one correspondence with edges
of the general map. Edges of the i-l.2.f.m strictly in $f_1$ (respectively $f_2$) correspond to edges 
of the general map in the first (respectively second) Vorono\"\i\  cell. As for edges on the loop separating $f_1$ and $f_2$, they come in three species: edges of type $m\, \rule[1.5pt]{8.pt}{1.pt}\, m$ correspond to 
map edges of type $(m-1)\, \rule[1.5pt]{8.pt}{1.pt}\, (m-1)$ shared by the two cells and receive the weight $\sqrt{g\, h}$ accordingly; edges of type $m\, \rule[1.5pt]{8.pt}{1.pt}\, (m+1)$ (when oriented with $f_1$ on their left) correspond to edges of the general map of type $m\, \rule[1.5pt]{8.pt}{1.pt}\, (m-1)$ in the first cell and 
edges of type $(m+1)\, \rule[1.5pt]{8.pt}{1.pt}\, m$ correspond to edges of the general map of type $m\, \rule[1.5pt]{8.pt}{1.pt}\, (m-1)$ in the second cell. We are thus lead to assign the weight $g$ to loop 
edges of 
the second species and $h$ to loop edges of the third species but, since there is clearly the same number of edges of the two types in a closed loop, we way equivalently assign the weight $\sqrt{g\, h}$ to all of them.

Again, writing $\delta(v_1,v_2)=2s$ for general maps enumerated by $F^{\rm{even}}(g,h)$ and $\delta(v_1,v_2)=2s-1$ for general maps enumerated by $F^{\rm{odd}}(g,h)$, with $s\geq 1$, we may decide to shift all labels by $-s$ in the associated  i-l.2.f.m. With these shifted labels, the planar i-l.2.f.m enumerated by $F^{\rm{even}}(g,h)$ may alternatively be characterized by the same rules
$\langle \hbox{c}_2\rangle$-$\langle \hbox{c}_4\rangle$ as before but with $\langle \hbox{c}_1\rangle$ replaced by the slightly more restrictive rule:.
\begin{enumerate}[$\langle \hbox{c}_1\rangle$-even:]
\item{The minimal label for the set of loop vertices is $0$ and none of the loop edges has labels $0\, \rule[1.5pt]{8.pt}{1.pt}\, 0$. The edges of the loop receive a weight $\sqrt{g\, h}$.}
\end{enumerate}
Similarly, for planar i-l.2.f.m enumerated by $F^{\rm{odd}}(g,h)$, $\langle \hbox{c}_1\rangle$ is replaced by the rule:
\begin{enumerate}[$\langle \hbox{c}_1\rangle$-odd:]
\item{The minimal label for the set of loop vertices is $0$ and at least one loop edge has labels $0\, \rule[1.5pt]{8.pt}{1.pt}\, 0$. The edges of the loop receive a weight $\sqrt{g\, h}$.}
\end{enumerate}

The conditions $\langle \hbox{c}_1\rangle$-even and $\langle \hbox{c}_1\rangle$-odd are clearly complementary among i-l.2.f.m satisfying the condition $\langle \hbox{c}_1\rangle$.
We immediately deduce that
\begin{equation*}
F^{\rm{even}}(g,h)+F^{\rm{odd}}(g,h)=F(g,h)
\end{equation*}   
so that we may interpret $F(g,h)$ as the generating function for bi-pointed general planar maps with two marked vertices at \emph{arbitrary distance} from each other,
with a weight $g$ per edge in the first Vorono\"\i\ cell, $h$ per edge in the second cell, and $\sqrt{g\, h}$ per edge shared by both cells. 
As a direct consequence, among bi-pointed general planar maps of fixed area $N$, with their two marked vertices at arbitrary distance, the law for the ratio $\phi=n/N$ of the area $n$ 
of one of the two Vorono\"\i\ cells by the total area $N$ is again, for large $N$, uniform between $0$ and $1$.

If we wish to control the parity of $\delta(v_1,v_2)$, we have to take into account the new constraints on loop edges.
We invite the reader to look at \cite{FG14b} for a detailed discussion on how to incorporate these constraints.
For $\delta(v_1,v_2)$ even, the generating function $F^{\rm{even}}(g,h)$ may be written as 
\begin{equation*}
F^{\rm{even}}(g,h)=\sum_{s\geq 1} \Delta_s\Delta_t \log(N_{s,t}(g,h))\Big\vert_{t=s} = \sum_{s\geq 1} \log\left(\frac{N_{s,s}(g,h)N_{s-1,s-1}(g,h)}{N_{s-1,s}(g,h)N_{s,s-1}(g,h)}\right)
\end{equation*}
where
\begin{equation}
N_{s,t}(g,h)=\frac{X_{s,t}(g,h)}{1+\sqrt{g\, h}\, R_s(g)\, R_t(h)\, X_{s,t}(g,h)}
\label{eq:NX}
\end{equation}
enumerates l.c with none of their spine edges having labels $0\, \rule[1.5pt]{8.pt}{1.pt}\, 0$.

For $\delta(v_1,v_2)$ odd, the generating function $F^{\rm{odd}}(g,h)$ reads (see again \cite{FG14b})
\begin{equation*}
F^{\rm{odd}}(g,h)=\sum_{s\geq 1} \Delta_s\Delta_t \log\left(\frac{X_{s,t}(g,h)}{N_{s,t}(g,h)}\right)\Big\vert_{t=s} =F(g,h)-F^{\rm{even}}(g,h)
\end{equation*}
as it should.

We may estimate the singularity $F^{\rm{even}}(g,h)$ from the scaling function 
associated with $N_{s,t}(g,h)$ and from its value at $g=h=1/12$. 
It is easily checked from its expression \eqref{eq:NX} that
\begin{equation*}
N_{\left\lfloor S/\epsilon\right\rfloor,\left\lfloor T/\epsilon\right\rfloor}(g,h)  =\frac{3}{2}+\frac{1}{4} x(S,T,a,b)\ \epsilon+O(\epsilon^2)
\end{equation*}
and, by the same arguments as for quadrangulations,
\begin{equation*}
\begin{split}
F^{\rm{even}}(g,h)
& =\epsilon^2 \int_{\epsilon}^{\infty} dS\,  \frac{1}{6} \partial_S\partial_T x(S,T,a,b)\Big\vert_{T=S} \\
&+\sum_{s=1}^{\left\lfloor S_0/\epsilon\right\rfloor-1} \log \left(\frac{(2s+1)^3 (2 s+3)}{(2 s+2)^3 \, 2s}\right)-\epsilon^2 \int_\epsilon^{S_0}dS\, \frac{1}{4S^3}+O(\epsilon^3)\ .
\\
\end{split}
\end{equation*}
This yields the expansion:
\begin{equation*}
\begin{split}
F^{\rm{even}}(g,h)
= \log \left(\frac{32}{3\pi^2}\right) -\frac{1}{36}\, \frac{(a^6-b^6)}{(a^4-b^4)}\, \epsilon^2+O(\epsilon ^3)\ .\\
\end{split}
\end{equation*}
with, as expected, the same singularity as $F(g,h)$ up to a factor $1/2$ since the number of bi-pointed general maps with $\delta(v_1,v_2)$ even is (asymptotically) half
the number of bi-pointed quadrangulations with $d(v_1,v_2)$ even.
Again, for the restricted ensemble of bi-pointed general planar map whose marked vertices are at even distance from each other, the law for the ratio $\phi=n/N$ of the area $n$ of one of the two Vorono\"\i\ cells by the total area $N$ is, for large $N$, uniform between $0$ and $1$. The same is obviously true if we condition the distance to be odd since
\begin{equation*}
\begin{split}
F^{\rm{odd}}(g,h)
= \log \left(\frac{\pi^2}{8}\right) -\frac{1}{36}\, \frac{(a^6-b^6)}{(a^4-b^4)}\, \epsilon^2+O(\epsilon ^3)\ .\\
\end{split}
\end{equation*}

\section{Conclusion}
\label{sec:conclusion}
In this paper, we computed the law for the ratio $\phi=n/N$ of the area $n$ ($=$ number of faces) of one of the two Vorono\"\i\ cells by the total area $N$ for random planar quadrangulations
with a large area $N$ and two randomly chosen marked distinct vertices at even distance from each other. We found that this law is uniform between $0$ and $1$,
which corroborates Chapuy's conjecture. We then extended this result to the law for the ratio $\phi=n/N$ of the area $n$ ($=$ number of edges) of one of the two Vorono\"\i\ cells by the total area $N$ for random general planar maps with a large area $N$ and two randomly chosen marked distinct vertices at arbitrary distance from each other. 
We again found that this law is uniform between $0$ and $1$.

Our calculation is based on an estimation of the singularity of the appropriate generating function keeping a control on the area of the Vorono\"\i\ cells,
itself based on an estimation of the singularity of some particular generating function $X_{s,t}(g,h)$ for labelled chains. Clearly, a challenging problem
would be to find an exact expression for $X_{s,t}(g,h)$ as it would certainly greatly simplify our derivation. 

Chapuy's conjecture extends to an arbitrary number of Vorono\"\i\ cells in a map of arbitrary fixed genus. It seems possible to test it by our method for
some slightly more involved cases than the one discussed here, say with three Vorono\"\i\ cells in the planar case or for two Vorono\"\i\ cells in maps with genus $1$. An important 
step toward this calculation would be to estimate the singularity of yet another generating function, $Y_{s,t,u}(g,h,k)$ enumerating labelled trees
with three non-aligned marked vertices and a number of label constraints\footnote{See \cite{BG08} for a precise list of label constraints.
There $Y_{s,t,u}(g,h,k)$ is defined when $g=h=k$ but the label constraints are independent of the weights.} involving subtrees divided into three
subsets with edge weights $g$, $h$, and $k$ respectively. Indeed, applying the Miermont bijection to maps with more points or for higher genus creates
labelled maps whose ``skeleton" (i.e.\ the frontier between faces)  is no longer a single loop but has branching points enumerated by $Y_{s,t,u}(g,h,k)$.
This study will definitely require more efforts.     

Finally, in view of the simplicity of the conjectured law, one may want to find a general argument which makes no use of any precise enumeration result
but relies only on bijective constructions and/or symmetry considerations.

\section*{Acknowledgements} 
I thank Guillaume Chapuy for bringing to my attention his nice conjecture and Timothy Budd for clarifying discussions. I also acknowledge the support of the grant ANR-14-CE25-0014 (ANR GRAAL). 

\newpage
\appendix
\section{Expression for the scaling function $x(S,T,a,b)$} 
The scaling function $x(S,T,a,b)$, determined by the partial differential equation 
\begin{equation*}
2 \big(x(S,T,a,b)\big)^2+6 \big(\partial_Sx(S,T,a,b)+\partial_Tx(S,T,a,b)\big)+27 \big(r(S,a)+r(T,b)\big)=0
\end{equation*}
(with $r(S,a)$ as in \eqref{eq:expRg}) and by the small $S$ and $T$ behavior \eqref{eq:smalST} is given by
\begin{equation*}
x(S,T,a,b)=-3\sqrt{\frac{a^2+b^2}{2}}-\frac{\mathfrak{N}(e^{-a\, S},e^{-b\, T},a,b)}{(1-e^{-a\, S} ) (1-e^{-b\, T}) \mathfrak{D}(e^{-a\, S},e^{-b\, T},a,b)}\ ,
\end{equation*}
where the polynomials
\begin{equation*}
\mathfrak{N}(\sigma,\tau,a,b)=\sum_{i=0}^{3}\sum_{j=0}^3 n_{i,j}\, \sigma^i\tau^j  \quad \hbox{and}\quad
\mathfrak{D}(\sigma,\tau,a,b) =\sum_{i=0}^{2}\sum_{j=0}^2 d_{i,j}\, \sigma^i\tau^j 
\end{equation*}
have the following coefficients $n_{i,j}\equiv n_{i,j}(a,b)$ and $d_{i,j}\equiv d_{i,j}(a,b)$:
writing for convenience these coefficients in the form
\begin{equation*}
n_{i,j}=n_{i,j}^{(0)}+\sqrt{\frac{a^2+b^2}{2}}n_{i,j}^{(1)}\ , \qquad
d_{i,j}=d_{i,j}^{(0)}+\sqrt{\frac{a^2+b^2}{2}}d_{i,j}^{(1)}\ ,
\end{equation*}
we have
\begin{equation*}
\begin{split}
n_{0,0}^{(0)}&= n_{0,3}^{(0)} = n_{3,0}^{(0)} = n_{3,3}^{(0)}= 0\\ 
n_{0,1}^{(0)}&= -\frac{18 b^3}{2 a^2+b^2} \qquad n_{0,2}^{(0)}= \frac{18 b^3 \left(5 a^2 +7 b^2\right)}{(a-b) (a+b) \left(2 a^2+b^2\right)}\\ 
n_{1,0}^{(0)}&= -\frac{18 a^3}{a^2+2 b^2} \qquad n_{2,0}^{(0)}= -\frac{18 a^3 \left(7 a^2+5 b^2\right)}{(a-b) (a+b) \left(a^2+2 b^2\right)}\\ 
n_{1,1}^{(0)}&= \frac{54 \left(2 a^7+17 a^5 b^2+17 a^4 b^3+17 a^3 b^4+17 a^2 b^5+2 b^7\right)}{(a-b)^2 \left(2 a^2+b^2\right) \left(a^2+2 b^2\right)}\\ 
n_{1,2}^{(0)}&=-\frac{54 \left(2 a^7+8 a^6 b+27 a^5 b^2+47 a^4 b^3+47 a^3 b^4+51 a^2 b^5+20 a b^6+14 b^7\right)}{(a-b) (a+b) \left(2 a^2+b^2\right) \left(a^2+2 b^2\right)}\\ 
n_{1,3}^{(0)}&= \frac{18 a^2 \left(2 a^5+12 a^4 b+17 a^3 b^2+36 a^2 b^3+17 a b^4+24 b^5\right)}{(a-b) (a+b) \left(2 a^2+b^2\right) \left(a^2+2 b^2\right)}\\ 
n_{2,1}^{(0)}&= \frac{54 \left(14 a^7+20 a^6 b+51 a^5 b^2+47 a^4 b^3+47 a^3 b^4+27 a^2 b^5+8 a b^6+2 b^7\right)}{(a-b) (a+b) \left(2 a^2+b^2\right) \left(a^2+2 b^2\right)}\\ 
n_{2,2}^{(0)}&= -\frac{54 \left(14 a^7+12 a^6 b+41 a^5 b^2+41 a^4 b^3+41 a^3 b^4+41 a^2 b^5+12 a b^6+14 b^7\right)}{(a-b)^2 \left(2 a^2+b^2\right) \left(a^2+2 b^2\right)}\\ 
n_{2,3}^{(0)}&= \frac{18 a^2 \left(14 a^5+32 a^4 b+51 a^3 b^2+58 a^2 b^3+37 a b^4+24 b^5\right)}{(a-b)^2 \left(2 a^2+b^2\right) \left(a^2+2 b^2\right)}\\ 
n_{3,1}^{(0)}&= -\frac{18 b^2 \left(24 a^5+17 a^4 b+36 a^3 b^2+17 a^2 b^3+12 a b^4+2 b^5\right)}{(a-b) (a+b) \left(2 a^2+b^2\right) \left(a^2+2 b^2\right)}\\ 
n_{3,2}^{(0)}&= \frac{18 b^2 \left(24 a^5+37 a^4 b+58 a^3 b^2+51 a^2 b^3+32 a b^4+14 b^5\right)}{(a-b)^2 \left(2 a^2+b^2\right) \left(a^2+2 b^2\right)}\\ 
\end{split}
\end{equation*}
and
\begin{equation*}
\begin{split}
n_{0,0}^{(1)}&= n_{0,3}^{(1)}= n_{3,0}^{(1)}= n_{3,3}^{(1)}= 0\\ 
n_{0,1}^{(1)}&= \frac{36 b^2}{2 a^2+b^2}\qquad n_{0,2}^{(1)}= -\frac{36 b^2 \left(a^2+5 b^2\right)}{(a-b) (a+b) \left(2 a^2+b^2\right)}\\ 
n_{1,0}^{(1)}&=\frac{36 a^2}{a^2+2 b^2} \qquad n_{2,0}^{(1)}= \frac{36 a^2 \left(5 a^2+b^2\right)}{(a-b) (a+b) \left(a^2+2 b^2\right)}\\ 
n_{1,1}^{(1)}&= -\frac{216 \left(a^6-a^5 b+8 a^4 b^2+2 a^3 b^3+8 a^2 b^4-a b^5+b^6\right)}{(a-b)^2 \left(2 a^2+b^2\right) \left(a^2+2 b^2\right)}\\ 
n_{1,2}^{(1)}&=\frac{216 \left(a^2+a b+b^2\right) \left(a^4+a^3 b+9 a^2 b^2+2 a b^3+5 b^4\right)}{(a-b) (a+b) \left(2 a^2+b^2\right) \left(a^2+2 b^2\right)}\\ 
n_{1,3}^{(1)}&= -\frac{36 a^2 \left(2 a^4+6 a^3 b+17 a^2 b^2+12 a b^3+17 b^4\right)}{(a-b) (a+b) \left(2 a^2+b^2\right) \left(a^2+2 b^2\right)}\\ 
n_{2,1}^{(1)}&= -\frac{216 \left(a^2+a b+b^2\right) \left(5 a^4+2 a^3 b+9 a^2 b^2+a b^3+b^4\right)}{(a-b) (a+b) \left(2 a^2+b^2\right) \left(a^2+2 b^2\right)}\\ 
n_{2,2}^{(1)}&=\frac{216 \left(5 a^6+4 a^5 b+13 a^4 b^2+10 a^3 b^3+13 a^2 b^4+4 a b^5+5 b^6\right)}{(a-b)^2 \left(2 a^2+b^2\right) \left(a^2+2 b^2\right)} \\ 
n_{2,3}^{(1)}&=-\frac{36 a^2 \left(10 a^4+22 a^3 b+33 a^2 b^2+26 a b^3+17 b^4\right)}{(a-b)^2 \left(2 a^2+b^2\right) \left(a^2+2 b^2\right)}\\ 
n_{3,1}^{(1)}&= \frac{36 b^2 \left(17 a^4+12 a^3 b+17 a^2 b^2+6 a b^3+2 b^4\right)}{(a-b) (a+b) \left(2 a^2+b^2\right) \left(a^2+2 b^2\right)}\\ 
n_{3,2}^{(1)}&= -\frac{36 b^2 \left(17 a^4+26 a^3 b+33 a^2 b^2+22 a b^3+10 b^4\right)}{(a-b)^2 \left(2 a^2+b^2\right) \left(a^2+2 b^2\right)}\ ,\\ 
\end{split}
\end{equation*}
while
\begin{equation*}
\begin{split}
d_{0,0}^{(0)}&=1\\ 
d_{0,1}^{(0)}&= -\frac{4 \left(a^2+2 b^2\right)}{2 a^2+b^2}\\ 
d_{0,2}^{(0)}&= \frac{2 a^4+17 a^2 b^2+17 b^4}{(a-b) (a+b) \left(2 a^2+b^2\right)}\\ 
d_{1,0}^{(0)}&= -\frac{4 \left(2 a^2+b^2\right)}{a^2+2 b^2}\\ 
d_{1,1}^{(0)}&=\frac{8 \left(4 a^2+a b+4 b^2\right) \left(a^4+7 a^2 b^2+b^4\right)}{(a-b)^2 \left(2 a^2+b^2\right) \left(a^2+2 b^2\right)}\\ 
d_{1,2}^{(0)}&=-\frac{4 \left(4 a^5+14 a^4 b+22 a^3 b^2+32 a^2 b^3+19 a b^4+17 b^5\right)}{(a-b) \left(2 a^2+b^2\right) \left(a^2+2 b^2\right)}\\ 
d_{2,0}^{(0)}&= -\frac{17 a^4+17 a^2 b^2+2 b^4}{(a-b) (a+b) \left(a^2+2 b^2\right)}\\ 
d_{2,1}^{(0)}&=\frac{4 \left(17 a^5+19 a^4 b+32 a^3 b^2+22 a^2 b^3+14 a b^4+4 b^5\right)}{(a-b) \left(2 a^2+b^2\right) \left(a^2+2 b^2\right)}\\ 
d_{2,2}^{(0)}&=-\frac{34 a^6+76 a^5 b+137 a^4 b^2+154 a^3 b^3+137 a^2 b^4+76 a b^5+34 b^6}{(a-b)^2 \left(2 a^2+b^2\right) \left(a^2+2 b^2\right)} \\ 
\end{split}
\end{equation*}
and 
\begin{equation*}
\begin{split}
d_{0,0}^{(1)}&=0\\ 
d_{0,1}^{(1)}&= \frac{12 b}{2 a^2+b^2}\\ 
d_{0,2}^{(1)}&=-\frac{12 b \left(a^2+2 b^2\right)}{(a-b) (a+b) \left(2 a^2+b^2\right)} \\ 
d_{1,0}^{(1)}&= \frac{12 a}{a^2+2 b^2}\\ 
d_{1,1}^{(1)}&=-\frac{48 \left(a^2+a b+b^2\right) \left(a^4+7 a^2 b^2+b^4\right)}{(a-b)^2 (a+b) \left(2 a^2+b^2\right) \left(a^2+2 b^2\right)}\\ 
d_{1,2}^{(1)}&=\frac{12 \left(2 a^4+6 a^3 b+11 a^2 b^2+9 a b^3+8 b^4\right)}{(a-b) \left(2 a^2+b^2\right) \left(a^2+2 b^2\right)}\\ 
d_{2,0}^{(1)}&= \frac{12 a \left(2 a^2+b^2\right)}{(a-b) (a+b) \left(a^2+2 b^2\right)}\\ 
d_{2,1}^{(1)}&=-\frac{12 \left(8 a^4+9 a^3 b+11 a^2 b^2+6 a b^3+2 b^4\right)}{(a-b) \left(2 a^2+b^2\right) \left(a^2+2 b^2\right)}\\
d_{2,2}^{(1)}&= \frac{12 (a+b) \left(a^2+a b+b^2\right) \left(4 a^2+a b+4 b^2\right)}{(a-b)^2 \left(2 a^2+b^2\right) \left(a^2+2 b^2\right)}\ .\\ 
\end{split}
\end{equation*}
It is easily verified that, for $b\to a$, $x(S,T,a,b)$ tends to $x(S,T,a)$ given by \eqref{eq:xa}, as expected.  

\section{Expression for the primitive $\mathfrak{K}(\sigma,a,b)$} 
Taking $\mathfrak{K}(\sigma,a,b)$ in the form 
\begin{equation*}
\mathfrak{K}(\sigma,a,b)=\frac{a\, b\, \sigma\, \tau\, \mathfrak{H}(\sigma,\tau,a,b)}{\left(\mathfrak{D}(\sigma,\tau,a,b)\right)^2}\Bigg\vert_{\tau=\sigma^{b/a}}
\end{equation*}
with the same function $\mathfrak{D}(\sigma,\tau,a,b)$ as in Appendix~A and where $\mathfrak{H}(\sigma,\tau,a,b)$ is a polynomial of the form
\begin{equation*}
\mathfrak{H}(\sigma,\tau,a,b)=\sum_{i=0}^{2}\sum_{j=0}^2 h_{i,j}\, \sigma^i\tau^j\ ,
\end{equation*}
the desired condition \eqref{eq:K} is fulfilled if   
\begin{equation*}
 \frac{1}{3}\partial_\sigma\partial_\tau \mathfrak{X}(\sigma,\tau,a,b)= \left\{(a+b)+a\, \sigma\, \partial_\sigma +b\, \tau\, \partial_\tau \right\}\frac{\mathfrak{H}(\sigma,\tau,a,b)}{\left(\mathfrak{D}(\sigma,\tau,a,b)\right)^2}\ .
 \end{equation*}
This fixes the coefficients $h_{i,j}\equiv h_{i,j}(a,b)$, namely:
\begin{equation*}
h_{i,j}=h_{i,j}^{(0)}+\sqrt{\frac{a^2+b^2}{2}}h_{i,j}^{(1)}
\end{equation*}
with
\begin{equation*}
\begin{split}
\hskip -1.2cm h_{0,1}^{(0)}&= 
 h_{1,0}^{(0)}= 
 h_{1,1}^{(0)}=
 h_{1,2}^{(0)}=
 h_{2,1}^{(0)}= 0\\ 
\hskip -1.2cm h_{0,0}^{(0)}&= -\frac{72 a^2 b^2 \left(4 a^2+a b+4 b^2\right)}{(a-b)^2 \left(2 a^2+b^2\right) \left(a^2+2 b^2\right)}\\ 
\hskip -1.2cm h_{0,2}^{(0)}&= \frac{72 a^2 b^2 \left(8 a^7+46 a^6 b+114 a^5 b^2+237 a^4 b^3+261 a^3 b^4+333 a^2 b^5+157 a b^6+140 b^7\right)}{(a-b)^3 (a+b)^2 \left(2 a^2+b^2\right)^2 \left(a^2+2 b^2\right)}\\  
\hskip -1.2cm h_{2,0}^{(0)}&= -\frac{72 a^2 b^2 \left(140 a^7+157 a^6 b+333 a^5 b^2+261 a^4 b^3+237 a^3 b^4+114 a^2 b^5+46 a b^6+8 b^7\right)}{(a-b)^3 (a+b)^2 \left(2 a^2+b^2\right) \left(a^2+2 b^2\right)^2}\\ 
\hskip -1.2cm h_{2,2}^{(0)}&= \frac{72 a^2 b^2 \left(4 a^2+a b+4 b^2\right) \left(70 a^6+148 a^5 b+281 a^4 b^2+298 a^3 b^3+281 a^2 b^4+148 a b^5+70 b^6\right)}{(a-b)^4 \left(2 a^2+b^2\right)^2 \left(a^2+2
   b^2\right)^2}\\ 
\end{split}
\end{equation*}
and
\begin{equation*}
\begin{split}
\hskip -1.2cm h_{0,1}^{(1)}&=
h_{1,0}^{(1)}=
h_{1,1}^{(1)}= 
h_{1,2}^{(1)}= 
h_{2,1}^{(1)}=0 \\ 
\hskip -1.2cm h_{0,0}^{(1)}&=\frac{432 a^2 b^2 \left(a^2+a b+b^2\right)}{(a-b)^2 (a+b) \left(2 a^2+b^2\right) \left(a^2+2 b^2\right)}\\  
\hskip -1.2cm h_{0,2}^{(1)}&=-\frac{432 a^2 b^2 \left(2 a^6+10 a^5 b+29 a^4 b^2+43 a^3 b^3+62 a^2 b^4+37 a b^5+33 b^6\right)}{(a-b)^3 (a+b)^2 \left(2 a^2+b^2\right)^2 \left(a^2+2 b^2\right)}\\ 
\hskip -1.2cm h_{2,0}^{(1)}&=\frac{432 a^2 b^2 \left(33 a^6+37 a^5 b+62 a^4 b^2+43 a^3 b^3+29 a^2 b^4+10 a b^5+2 b^6\right)}{(a-b)^3 (a+b)^2 \left(2 a^2+b^2\right) \left(a^2+2 b^2\right)^2}\\ 
\hskip -1.2cm h_{2,2}^{(1)}&=-\frac{1296 a^2 b^2 \left(a^2+a b+b^2\right) \left(22 a^6+52 a^5 b+89 a^4 b^2+106 a^3 b^3+89 a^2 b^4+52 a b^5+22 b^6\right)}{(a-b)^4 (a+b) \left(2 a^2+b^2\right)^2 \left(a^2+2
   b^2\right)^2}\ .\\ 
\end{split}
\end{equation*}
\section{Contour integral over $a$} 
\begin{figure}
\begin{center}
\includegraphics[width=6cm]{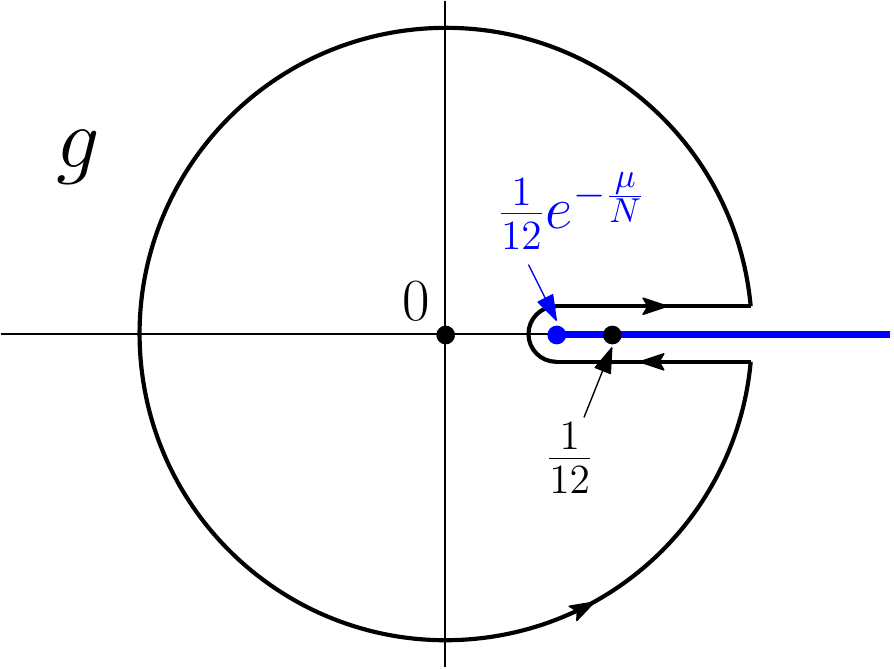}
\end{center}
\caption{Deformation of the contour for the integral over $g$.}
\label{fig:contourg}
\end{figure}
Given $\mu\geq 0$, the integral over $g$ 
\begin{equation*}
\frac{1}{2\rm{i}\pi} \oint\frac{dg}{g^{N+1}} F(g,g\, e^{\frac{\mu}{N}})
\end{equation*}
is on a contour around $0$. Here $F(g,g\, e^{\frac{\mu}{N}})$ has a singularity for real $g>(1/12) e^{-\frac{\mu}{N}}$ and the contour may be deformed as in Figure~\ref{fig:contourg}.
For large $N$, the dominant contribution comes from the vicinity of the cut and is captured by setting
\begin{figure}
\begin{center}
\includegraphics[width=11cm]{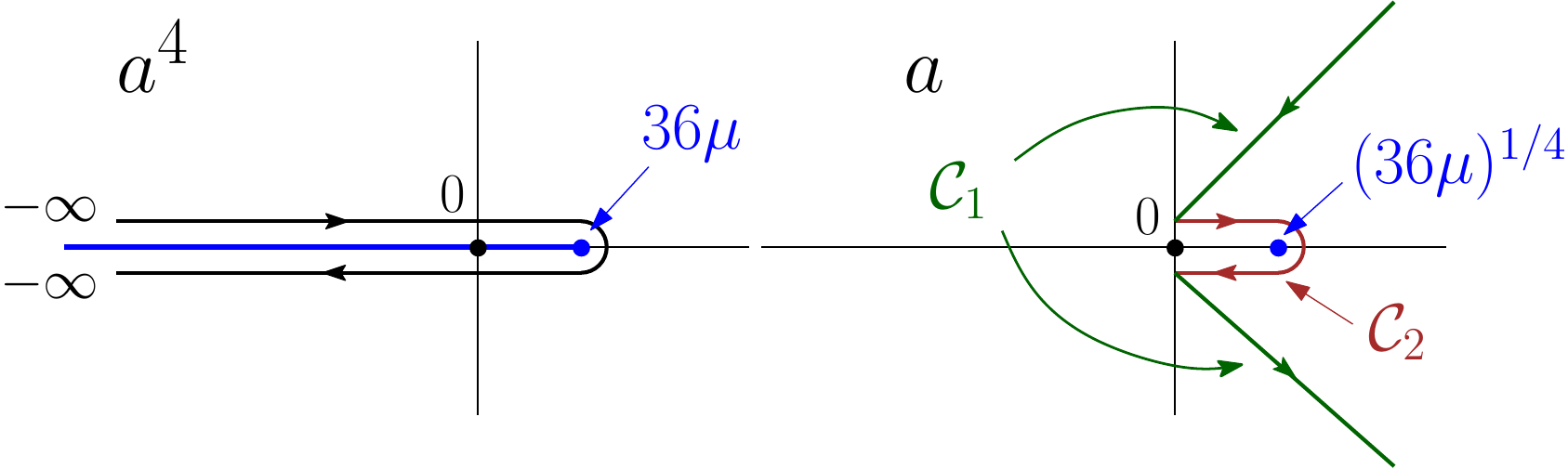}
\end{center}
\caption{The contour in the variable $a^4$ resulting from the large $N$ limit of the integral over $g$ along the contour of Figure~\ref{fig:contourg}. The resulting
countour $\mathcal{C}=\mathcal{C}_1\cap \mathcal{C}_2$ in the variable $a$.}
\label{fig:contoura}
\end{figure}
\begin{equation*}
g=\frac{1}{12}\left(1-\frac{a^4}{36}\, \frac{1}{N}\right)
\end{equation*}
where the variable $a^4$ varies along the cut from $-\infty$ to  $36 \mu$ back to $-\infty$. In other words, the contour $\mathcal{C}$ for the variable $a$ is that of Figure~\ref{fig:contoura},
made of two parts: a contour $\mathcal{C}_1$ made of two half straight lines at $\pm 45^\circ$ starting from the origin, and a contour $\mathcal{C}_2$ consisting of a back and 
forth excursion from $0$ to $(36\mu)^{1/4}$ back to $0$. In the variable $a$, the integral reads
\begin{equation*}
\frac{1}{2\rm{i}\pi} \int_\mathcal{C} da\,  \frac{-a^3}{9} \left\{ -\frac{1}{18} \frac{a^6-(a^4-36\mu)^{3/2})}{36\mu}\, e^{a^4/36} \right\}
\end{equation*}
Concerning the contour $\mathcal{C}_2$, the term $a^6$ has no cut hence contributes $0$ to the integral. As for the $(a^4-36\mu)^{3/2}$ term, setting 
$a=\sqrt{6}(\mu-t^2)^{1/4}$ with real $t$ from $\sqrt{\mu}$ to $0$ back to $\sqrt{\mu}$, we have
\begin{equation*}
\begin{split}
\frac{1}{2\rm{i}\pi} \int_{\mathcal{C}_2} da\,  \frac{-a^3}{9} \left\{\frac{1}{18} \frac{(a^4-36\mu)^{3/2})}{36\mu}\, e^{a^4/36} \right\}&=
\frac{1}{2\rm{i}\pi} \left\{\int_{\sqrt{\mu}}^0 dt\, 2t\, \frac{(-t^2)^{3/2}}{3\mu}\, e^{\mu-t^2} \right.
\\&\qquad \left.+\int_0^{\sqrt{\mu}} dt\, 2t\, \frac{(-t^2)^{3/2}}{3\mu}\, e^{\mu-t^2} \right\}\\
\end{split}
\end{equation*}
where $(-t^2)^{3/2}=-{\rm{i}}\,  t^3$ for the first integral and $(-t^2)^{3/2}={\rm{i}}\,  t^3$ for the second, so that the final contribution of the contour $\mathcal{C}_2$
is
\begin{equation*}
\frac{2}{3 \pi} \frac{e^{\mu}}{\mu} \int_0^{\sqrt{\mu}} dt\, t^4\, e^{-t^2}\ .
\end{equation*}
Let us now come to the integral over the contour $\mathcal{C}_1$.
The term $a^6$ now contributes to the integral: setting $a=\sqrt{6}\, e^{\pm {\rm i} \frac{\pi}{4}} \sqrt{t}$ with real $t$ from $+\infty$ to $0$ (respectively $0$ to $+\infty$), 
we get a contribution
\begin{equation*}
\begin{split}
\frac{1}{2\rm{i}\pi} \int_{\mathcal{C}_1} da\,  \frac{-a^3}{9} \left\{-\frac{1}{18}\ \frac{a^6}{36\mu}\, e^{a^4/36} \right\}&=
\frac{1}{2\rm{i}\pi} \left\{{\rm{i}}\int_{+\infty}^0 dt\, 2t\, \frac{t^3}{3\mu}\, e^{-t^2} \right.
\\& \qquad \qquad \left.-{\rm{i}} \int_0^{+\infty} dt\, 2t\, \frac{t^3}{3\mu}\, e^{-t^2} \right\}\\
&= -\frac{2}{3 \pi} \frac{1}{\mu} \int_0^{\infty} dt\, t^4\, e^{-t^2}\ .
\end{split}
\end{equation*}
Finally the $(a^4-36\mu)^{3/2}$ contribution is obtained by setting 
$a=\sqrt{6}\, e^{\pm {\rm i} \frac{\pi}{4}} (t^2-\mu)^{1/4}$ with real $t$ from $+\infty$ to $\sqrt{\mu}$ (respectively $\sqrt{\mu}$ to $+\infty$). We get
\begin{equation*}
\begin{split}
\frac{1}{2\rm{i}\pi} \int_{\mathcal{C}_1} da\,  \frac{-a^3}{9} \left\{\frac{1}{18} \frac{(a^4-36\mu)^{3/2})}{36\mu}\, e^{a^4/36} \right\}&=
\frac{1}{2\rm{i}\pi} \left\{\int_{+\infty}^{\sqrt{\mu}} dt\, 2t\, \frac{(-t^2)^{3/2}}{3\mu}\, e^{\mu-t^2} \right.
\\&\qquad \left.+\int_{\sqrt{\mu}}^{+\infty} dt\, 2t\,  \frac{(-t^2)^{3/2}}{3\mu}\, e^{\mu-t^2} \right\}\\
\end{split}
\end{equation*}
where again $(-t^2)^{3/2}=-{\rm{i}}\,  t^3$ for the first integral and $(-t^2)^{3/2}={\rm{i}}\,  t^3$ for the second, so that the final contribution reads
\begin{equation*}
\frac{2}{3 \pi} \frac{e^{\mu}}{\mu} \int_{\sqrt{\mu}}^\infty dt\, t^4\, e^{-t^2}\ .
\end{equation*}
Adding up all the contributions, we end up with the result:
\begin{equation*}
\begin{split}
\hskip -1.2cm \frac{2}{3 \pi} \frac{e^{\mu}}{\mu} \int_0^{\sqrt{\mu}} dt\, t^4\, e^{-t^2} -\frac{2}{3 \pi} \frac{1}{\mu} \int_0^{\infty} dt\, t^4\, e^{-t^2}
+\frac{2}{3 \pi} \frac{e^{\mu}}{\mu} \int_{\sqrt{\mu}}^\infty dt\, t^4\, e^{-t^2}& = \frac{e^\mu -1}{\mu} \frac{2}{3\pi}  \int_0^{\infty} dt\, t^4\, e^{-t^2}\\
& = \frac{e^\mu -1}{\mu} \times \frac{1}{4\, \sqrt{\pi}}\ .\\
\end{split}
\end{equation*}

\bibliographystyle{plain}
\bibliography{voronoi}

\end{document}